\input amstex
\input pictex
\UseAMSsymbols

%   \hoffset=0truecm \voffset=0truecm \hsize=12.7cm %\vsize=22.5cm
   \NoBlackBoxes
   \font\rmk=cmr8   \font\gross=cmbx10 scaled\magstep1   \font\abs=cmcsc10
   \font\itk=cmti8  \font\ttk=cmtt8
   \font\ttk=cmtt8
\def\t#1{\quad\text{#1}\quad}
\def\Hom{\operatorname{Hom}}
\def\Mod{\operatorname{Mod}}

\def\add{\operatorname{add}}
\def\Add{\operatorname{Add}}
\def\Ker{\operatorname{Ker}}
\def\Cok{\operatorname{Cok}}

\def\Ext{\operatorname{Ext}}
\def\pd{\operatorname{pd}}
\def\id{\operatorname{id}}
\def\mod{\operatorname{mod}}
\def\End{\operatorname{End}}
\def\Prod{\operatorname{Prod}}

   \def\Links{\abs I\. Reiten, C\.M\. Ringel}
   \def\Rechts{\abs Infinite dimensional representations of canonical algebras}
\headline{\ifnum\pageno=1\hfill %
    \else\ifodd\pageno \hfil\Rechts\hfil \else \hfil\Links\hfil \fi  \fi}
%\comment

\vglue 1cm \centerline{\gross Infinite Dimensional
Representations}
    \smallskip
\centerline{\gross of Canonical Algebras}
    \bigskip
\centerline{Idun Reiten and Claus Michael Ringel}

        \bigskip\medskip
\centerline{Dedicated to Vlastimil Dlab on the occasion of his
70th birthday.}

    \bigskip\bigskip

\plainfootnote{}{ {\rmk 2000 \itk Mathematics Subject
Classification. \rmk Primary
 16D70, %Modules: structure and classification
 16D90. %Module categories, Morita equivalence
Secondary
 16G20, %Representations of quivers
 16G60, %Representation type
 16G70. %Auslander-Reiten quivers
 }}
        \bigskip\bigskip
%===================================
{\baselineskip=9pt \narrower\narrower{\noindent \rmk ABSTRACT. The
aim of this paper is to extend the structure theory for infinitely
generated modules over tame hereditary algebras to the more
general case of modules over concealed canonical algebras. Using
tilting, we may assume that we deal with canonical algebras. The
investigation is centered around the generic and the Pr\"{u}fer
modules, and how other modules are determined by these
modules.}\par }

        \bigskip\bigskip
%===================================
{\bf Introduction.}
        \medskip
Let $\Lambda$ be a finite dimensional algebra over a field $k$.
Traditionally one mainly has considered the $\Lambda$-modules
which are finitely generated. An early exception were papers by
several authors dealing with modules over the Kronecker algebra.
This was generalized by the second author [R1] to the case of a
tame hereditary algebra $\Lambda$, an  investigation which was
based on the explicit knowledge of the finitely generated modules
as presented in his joint work [DR] with Dlab. As it turned out,
there are striking similarities between the category of all
$\Lambda$-modules and the category of all abelian groups (or the
category of all $R$-modules, where $R$ is a Dedekind ring with
infinitely many prime ideals). In particular, the so called
generic module and the Pr\"{u}fer modules play an important role,
as they correspond to the indecomposable injective $R$-modules.

The aim of the present paper is to show that the core results of
these old investigations only depend on the existence of a sincere
stable separating tubular family, and not at all on the
representation type of the algebra. Hence, in view of the
characterization due to Lenzing and de la Pe\~na in [LP], the
natural setting is the class of concealed canonical algebras,
which contains the class of tame hereditary algebras, but also
many others. An important special class is the better known class
of the canonical algebras, and actually, it is sufficient to deal
with this class (with the tubular family considered to be given by
the modules of defect zero), since it is easy to extend the
results to the general class of concealed canonical algebras via a
tilting procedure. Note that such a canonical algebra may be
domestic, or non-domestic tame, or wild, but always we will obtain
splitting results which are similar to those known for tame
hereditary algebras. In the special case of a canonical algebra
which is non-domestic tame (thus for all tubular algebras), there
are countably many tubular families: any such family gives rise to
corresponding split torsion pairs.

The key results are, as for tame hereditary algebras, centered
around the explicit description of some modules, which are defined
in a similar way as in the tame hereditary case: the generic
module and the Pr\"{u}fer modules. In some sense, all other
modules are determined from these, via maps between them. In order
to be more explicit, we need to introduce some notation and
terminology.

Given a ring $R$, we consider usually left $R$-modules and call
them just modules or also representations of $R$. The category of
all $R$-modules will be denoted by $\Mod R$, the full subcategory
of the finitely presented ones by $\mod R.$ For any class $\Cal X$
of $R$-modules, we denote by $\add\Cal X$ its additive closure: it
is the smallest full subcategory closed under isomorphisms, direct
summands and finite direct sums. Similarly, $\Add \Cal X$ is the
smallest full subcategory closed under isomorphisms, direct
summands and arbitrary direct sums, whereas $\Prod \Cal X$ is the
smallest full subcategory closed under isomorphisms, direct
summands and arbitrary products. Given $R$-modules $X,Y$, we
usually write $\Hom(X,Y)$ or $\Ext^1(X,Y)$ instead of
$\Hom_R(X,Y)$ or $\Ext^1_R(X,Y)$. When dealing with classes $\Cal
X, \Cal Y$ (or full subcategories) of $R$-modules, we write
$\Hom(\Cal X,\Cal Y) = 0$ in order to assert that $\Hom(X,Y) = 0$
for all $X\in \Cal X$ and $Y\in\Cal Y$, and similarly for
$\Ext^1.$ For any $R$-module $M$, we denote by $\pd M$ its
projective dimension and by $\id M$ its injective dimension. When
dealing with module classes (or full subcategories), two different
types of notations will be used: the module classes denoted by
script letters such as $\Cal X, \Cal C, \Cal Q$  (or also $\omega$
and $\omega_0$) will usually be closed under direct sums (often
even infinite direct sums); in contrast, when dealing with an
artin algebra, we will use small boldface letters such as $\bold
x, \bold p, \bold t$ in order to denote classes consisting only of
indecomposable modules of finite length.

Let $\Lambda$ be an artin algebra, and assume that there exist
classes $\bold p, \bold t, \bold q$ in $\mod \Lambda$ (a
``trisection'') with the following properties: $\bold t$ is a
sincere stable separating tubular family and it separates $\bold
p$ from $\bold q$ (see section 2). Note that an indecomposable
$\Lambda$-module of finite length belongs to $\bold p$ or $\bold
t$ if and only if it is cogenerated by $\bold t$. A crucial result
of this paper will be the following: {\it Any (not necessarily
finite dimensional) $\Lambda$-module $M$ has a direct sum
decomposition $M = M_0\oplus M_1$, where $M_0$ is cogenerated by
the direct limit closure $\Cal T$ of $\bold t$ and $M_1$ is
generated by $\bold t$;} in addition, we can assume that
$\Hom(M_1,\Cal T) = 0,$ and then $\Hom(M_1,M_0) = 0$. We consider
the class $\Cal C$ of modules cogenerated by $\Cal T$, thus a
finite length module belongs to $\Cal C$ if and only if it is
cogenerated by $\bold t$. It follows that $\Cal C$ is the
torsionfree class of a split torsion pair in $\Mod \Lambda$. We
will investigate in detail {\bf all} the torsion pairs in
$\Mod\Lambda$ with the property that a finite length module is
torsionfree if and only if it is cogenerated by $\bold t$. As we
will see, all these torsion pairs split. Now $\Cal C$ is the
largest possible torsionfree class of this kind. Also the largest
possible torsion class $\Cal D$ of such a torsion pair can be
described easily: it is the class of all modules $M$ with
$\Hom(M,\bold t) = 0.$  The category $\omega=\Cal C \cap \Cal D$
turns out to be of central importance. The main results of the
paper can be expressed in terms of these categories $\Cal C, \Cal
D$ and $\omega$. The objects in $\omega$ can be completely
classified: any object in $\omega$ is a direct sum of copies of
the generic module $G$ and of Pr\"{u}fer modules. The class $\Cal
C$ is determined by $\omega$ as $\{C\mid\Ext^1(C,\omega)=0\}$, and
$\Cal D$ is determined by $\omega$ as
$\{D\mid\Ext^1(\omega,D)=0\}$. Further there are exact sequences
$0\rightarrow C\rightarrow V\rightarrow V'\rightarrow 0$ with
$V\in$ $\omega$ and $V'$ a direct sum of Pr\"{u}fer modules, for
$C$ in $\Cal C$, and $0\rightarrow V'\rightarrow V\rightarrow D\to
0$, with $V'\in \Add G$ and $V\in \omega$, for $D$ in $\Cal D$. As
a consequence, the modules in $\Cal C$ can be characterized as the
kernels of maps in $\omega$, and similarly, the modules in $\Cal
D$ can be characterized as the cokernels of maps in $\omega.$ Thus
any $\Lambda$-module $M$ is obtained as a direct sum $M =
M_0\oplus M_1$, where $M_0$ is the kernel and $M_1$ the cokernel
of suitable maps in $\omega$: in this way, the category
$\Mod\Lambda$ can be completely described in terms of $\omega$.

When dealing with finite dimensional algebras, one may argue that
it is the category of finite dimensional representations which is
the primary object of interest. However, the relevance of infinite
dimensional representations has been stressed at various occasions
[R1,R8] and here we encounter again such a situation: it is the
subcategory $\omega$ which plays the decisive role when studying
the cut between $\bold t$ and $\bold q$ in $\mod \Lambda$, and as
we have noted, $\omega$ does not contain a single non-zero
finite-dimensional representation.  We will denote by $W$ the
direct sum of all the indecomposables in $\omega$, one from each
isomorphism class. This module $W$ allows to reconstruct $\omega$
(as $\Add W$), thus the whole category $\Mod\Lambda$. Clearly, $W$
is a very valuable module! This can be phrased quite well in terms
of tilting and cotilting theory. We will use the denomination
inf-tilting and inf-cotilting when we deal with the general
concepts without the restriction of dealing with finite
dimensional modules, see section 11. Our results show that $W$ is
both an inf-tilting module of projective dimension one and an
inf-cotilting module of injective dimension one (see [BS] for a
different and independent approach to this for cotilting modules).
It is also possible to perform tilting with respect to torsion
pairs as in [HRS] to construct new hereditary categories where the
objects in $\omega$ become enough projective or enough injective
objects.

If we consider the special case of a tame hereditary algebra, most
of the results presented here have been established in [R1], but
for Proposition 4 (the classification of torsion pairs) we should
refer to unpublished information by Assem and Kerner. It should be
noted that Theorem 5 (the existence of the right
$\omega$-approximations) seems to be new even in this case. The
proof is inspired by [AB].

We will follow quite closely the presentation given in [R1], using
only the structure theory for finite dimensional representations,
and not taking into account the large amount of information on
infinite dimensional representations obtained in the meantime by
various authors. In particular, we will construct the relevant
''generic'' module $G$ from scratch. At the end we indicate
 a different approach using the available results. The reader should not
mind that the text itself avoids all more sophisticated
considerations, but this stubborn approach should make it quite
transparent to trace in which way the structure of the category of
finite dimensional representations determines that of all the
representations.

The paper is organized as follows. In section 1 we give a
criterion for a torsion pair to be split. In section 2 we recall
basic properties of the central algebras in this paper; the
canonical and concealed canonical algebras. In the next six
sections (3-8) we deal with a canonical algebra $\Lambda$ and the
canonical trisection $(\bold p, \bold t,\bold q)$ of $\mod
\Lambda$. We investigate the two extremal torsion pairs of $\Mod
\Lambda$ mentioned above in section 3, and give the structure of
the Pr\"{u}fer modules. The left $\omega$-approximation sequence
is established in section 4, and the basic splitting result
$\Ext^1(\Cal C, \Cal D)=0$ is given in section 5. The structure of
$\omega$ is presented in section 6, and the existence of the right
$\omega$-approximation sequences is deduced in section 7. The
structure of $\omega$ is investigated more closely in section 8.
In section 9 we outline that all these considerations are valid
for any sincere stable separating tubular family, thus for any
concealed canonical algebra. Of course, we use tilting functors in
order to relate an arbitrary sincere stable separating tubular
family with the canonical trisection of a canonical algebra.
Connections with tilting theory are discussed in sections 10 and
11. In section 12 we provide further comments and indicate another
approach to the results in this paper. A tubular algebra has a lot
of sincere stable separating tubular families and as we will see
in section 13, our considerations allow to attach a non-negative
real number as a ``slope'' to any indecomposable infinite
dimensional module.
        \medskip
The investigations presented here have for the most part been
completed during a stay of the first author at Bielefeld in 1998
and she would like to thank the second author for his hospitality;
unfortunately, the write-up of the results has been delayed for
quite a while.
         \bigskip\bigskip
%=================================================
{\bf 1\. Torsion pairs.}
    \medskip
The investigations presented in this paper are centered around
various torsion pairs (or, as they are sometimes called, torsion
``theories''). We are going to recall the relevant definitions and
main properties, and we provide a general method for producing
split torsion pairs.
        \medskip
Let $R$ be a ring. For any class $\Cal Z$ of $R$-modules, we
denote by $l(\Cal Z)$ the class of all $R$-modules $M$ with
$\Hom(M,\Cal Z) = 0,$ and similarly, $r(\Cal Z)$ is the class of
all $R$-modules $M$ with $\Hom(\Cal Z,M) = 0$ (let us stress that
our notation $r(-)$ and $l(-)$ always refers to the complete
category $\Mod R$ as ambient category, the only exception being
section 10 where the ambient category is an arbitrary abelian
category).
        \medskip
{\bf Lemma 1.} {\it Let $\Cal F,\Cal G$ be classes of $R$-modules.
The following conditions are equivalent:
\item{\rm (i)} $l(\Cal F) = \Cal G$ and $r(\Cal G) = \Cal F.$
\item{\rm (ii)} $\Hom(\Cal G,\Cal F) = 0$ and any module $M$ has a submodule $M' \in \Cal G$
    such that $M/M' \in \Cal F$.}
        \medskip
If these conditions are satisfied, the pair $(\Cal F,\Cal G)$ is
said to be a {\it torsion pair} with {\it torsionfree class} $\Cal
F$, and  {\it torsion class} $\Cal G$. The modules in $\Cal F$ are
called the {\it torsionfree}, those in $\Cal G$ the {\it torsion}
modules. It is straightforward to see that the submodule $M'$
given in (ii) is uniquely determined by $M$ (provided the torsion
pair $(\Cal F,\Cal G)$ is fixed).
        \medskip
Proof of the equivalence. (i) $\implies$ (ii): Only the last
assertion needs a proof. Thus, let $M$ be an arbitrary $R$-module.
Let $M'$ be the sum of images of maps from a module in $\Cal G$ to $M$.
Since $\Cal G = l(\Cal F)$, $\Cal G$ is closed under factors and arbitrary sums,
so that $M'$ is in $\Cal G$. Since $\Cal G$ is also closed under extensions,
we see that $\Hom(\Cal G, M/M') = 0$, so that $M/M'$ is in $\Cal F$.

(ii) $\implies$ (i). We show that $l(\Cal F) \subseteq \Cal G$.
Let $N$ belong to $l(\Cal F)$. According to (ii)  there exists a
submodule $N'$ of $N$ which belongs to $\Cal G$ such that $N/N'$
belongs to $\Cal F$. But the assumption that $N \in l(\Cal F)$
implies that the projection map $N \to N/N'$ is the zero map, thus
$N/N' = 0$ and therefore $N = N' \in \Cal G$. Similarly, one shows
that $r(\Cal G) \subseteq \Cal F.$
        \medskip
Some readers may wonder about the not quite usual sequence of
naming the {\it torsionfree class} $\Cal F$ first and the {\it
torsion class} $\Cal G$ second
--- this corresponds to the vision of drawing arrows and thus non-trivial maps from left to right
(whenever possible): there usually will be many non-zero maps from
the objects in $\Cal F$ to the objects in $\Cal G$ (but, by
definition, none in the other direction), thus $\Cal F$ may be
considered as ``situated to the left'' of $\Cal G$.
        \medskip
The torsion pair $(\Cal F,\Cal G)$ is said to be {\it split}
provided $\Ext^1(\Cal F,\Cal G) = 0$, or, equivalently, provided
every module is the direct sum of a module in $\Cal F$ and a
module in $\Cal G$.
        \bigskip
Any class $\Cal Z$ of $R$-modules determines two torsion pairs,
namely $$
  (r(\Cal Z),lr(\Cal Z))\t{and} (rl(\Cal Z),l(\Cal Z)).
$$ Clearly, $lr(\Cal Z)$ is the smallest possible torsion class
containing $\Cal Z$, whereas $rl(\Cal Z)$ is the smallest possible
torsionfree class containing $\Cal Z.$
        \medskip
{\bf Lemma 2.} {\it Let $\Cal Z$ be any class of $R$-modules. Then
an $R$-module $M$ belongs to $lr(\Cal Z)$ if and only if the only
submodule $U$ of $M$ with $M/U\in r(\Cal Z)$ is $U = M$.
Similarly, an $R$-module $M$ belongs to $rl(\Cal Z)$ if and only
if the only submodule $U$ of $M$ with $U\in l(\Cal Z)$ is $U =
0$.}
        \medskip
Proof: If $M$ belongs to $lr(\Cal Z)$ and $U$ is a submodule of
$M$ with $M/U\in r(\Cal Z)$, then the projection map $M \to M/U$
has to be the zero map, thus $U = M$. Conversely, assume that $M$
is an $R$-module such that the only submodule $U$ with $M/U \in
r(\Cal Z)$ is $U = M$. Since $(r(\Cal Z),lr(\Cal Z))$ is a torsion
pair, the module $M$ has a submodule $M'$ which belongs to
$lr(\Cal Z)$ such that $M/M'$ belongs to $r(\Cal Z)$. Since $M'$
is a submodule of $M$ with $M/M' \in r(\Cal Z)$, we know by
assumption that $M'= M$. But this shows that $M \in lr(\Cal Z).$
This proves the first equivalence. The second equivalence is shown
in the same way.
        \medskip
Given a class $\Cal Z$ of $R$-modules, we denote by $g(\Cal Z)$
the class of all $R$-modules generated by $\Cal Z$ (these are just
the factors of direct sums of modules in $\Cal Z$, and by $c(\Cal
Z)$ those cogenerated by $\Cal Z$ (these are the submodules of
products of modules in $\Cal Z$). The following inclusions are
trivial: $$
 g(\Cal Z) \subseteq lr(\Cal Z) \t{and} c(\Cal Z) \subseteq rl(\Cal Z).
$$
        \medskip
{\bf Lemma 3.} {\it Let $\Cal Z$ be a class of $R$-modules. Then
$g(\Cal Z) = lr(\Cal Z)$ if and only if $g(\Cal Z)$ is closed
under extensions. Similarly, $c(\Cal Z)= rl(\Cal Z)$ if and only
if $c(\Cal Z)$ is closed under extensions.}
        \medskip
Proof: We show the first assertion  (the second assertion is shown
in the same way). Note that $lr(\Cal Z)$ is closed under
extensions, thus the equality $g(\Cal Z) = lr(\Cal Z)$ implies
that $g(\Cal Z)$ is closed under extensions. Conversely, assume
that $g(\Cal Z)$ is closed under extensions and let $M \in lr(\Cal
Z)$. We have to show that $M$ belongs to $g(\Cal Z)$. Let $M'$ be
the sum of all images of maps $Z \to M$ with $Z \in \Cal Z$, thus
$M'$ is the maximal submodule of $M$ generated by $\Cal Z$. We
claim that $M/M'$ belongs to $r(\Cal Z).$ Namely, given a map
$f\:Z \to M/M'$ with $Z\in \Cal Z$, let $M''/M'$ be its image,
where $M' \subseteq M'' \subseteq M$. Now $M'$ and $M''/M'$ are
generated by $\Cal Z$, thus, by assumption also $M''$ is generated
by $\Cal Z$. But this means that $M'' \subseteq M'$ and therefore
$f = 0.$ Since $M \in lr(\Cal Z)$ and $M/M' \in r(\Cal Z)$, the
projection map $M \to M/M'$ is the zero map, thus $M = M'$ (of
course, one also may refer to Lemma 2). This shows that $M$
belongs to $g(\Cal Z)$.
        \bigskip
It will be useful to know conditions so that a subcategory of the
form $g(\Cal Z)$ is closed under extensions. From now on, we
restrict to the case when $R = \Lambda$ is an artin algebra and we
will denote the Auslander-Reiten translation in $\mod \Lambda$ by
$\tau$. Let us consider the case when $\Cal Z = \bold z$ is a
class of modules of finite length.
        \medskip
{\bf Lemma 4.} {\it Let $\Lambda$ be an artin algebra and $\bold
z$ a class of $\Lambda$-modules of finite length. Assume that
$\add \bold z$ is closed under extensions. If either $\add\bold z$
is also closed under factor modules or if $\pd Z \le 1$ for all
$Z\in \bold z$, then $g(\bold z)$ is closed under extensions.}
        \medskip
Proof. We first show the following: Under either assumption, given
a finite length module $Y$ and a submodule $X$ of $Y$ such that
both $X$ and $Y/X$ are generated by $\bold z$, then also $Y$ is
generated by $\bold z$. Namely, if we assume that $\add\bold z$ is
closed under factor modules, then both $X$ and $Y/X$ belong to
$\add\bold z$, since they are factor modules of modules in
$\add\bold z$. Thus also $Y$ belongs to $\add\bold z$, since we
assume that $\add\bold z$ is closed under extensions. Next, assume
that $\pd Z \le 1$ for all $Z \in \bold z$. There are surjective
maps $\pi\:Z \to Y/X$ and $\pi'\:Z' \to X$ where $Z,Z'$ belong to
$\add \bold z.$ Starting from the exact sequence $0 \to X \to Y
\to Y/X \to 0,$ we can form the induced exact sequence with
respect to $\pi$. Using now that $\pd Z \le 1$, and that $\pi'$ is
an epimorphism, we obtain a commutative diagram with exact rows of
the following shape: $$ \CD
 0 @>>> X @>>> Y @>>> Y/X @>>> 0 \cr
 @.     @|    @AAf A     @AA\pi A  \cr
 0 @>>> X @>>> Y'@>>> Z    @>>> 0 \cr
 @.     @AA\pi' A  @AAf' A    @|  \cr
 0 @>>> Z' @>>> Y''@>>> Z    @>>> 0
\endCD
$$ On one hand, the map $ff'$ is surjective, on the other hand,
$Y''$ belongs to $\add\bold z$, since $\add\bold z$ is closed
under extensions. This shows that $Y$ is generated by $\bold z.$

Now consider the general case of an arbitrary $\Lambda$-module $Y$
and a submodule $X$ of $Y$ such that both $X$ and $Y/X$ are
generated by $\bold z$. We have to show that $Y$ is generated by
$\bold z$. Write $Y = \sum_i Y_i$, where $X \subseteq Y_i$ and
$Y_i/X$ is isomorphic to a factor module of some module in $\bold
z$. It is sufficient to show that all the $Y_i$ belong to $g(\bold
z)$. Thus, without loss of generality, we may assume that $Y/X$ is
of finite length. Since $Y/X$ is of finite length, there is a
finite length submodule $Y'$ of $Y$ with $Y = X+Y'$. Now, $X$ is
the filtered union of submodules $X_i$ of finite length generated
by $\bold z$, thus there is some $i$ with $X\cap Y' = X_i \cap
Y'.$ Thus $(X_i+Y')/X_i \simeq Y'/(X_i\cap Y') = Y'/(X\cap Y')
\simeq (X+Y')/X = Y/X$. This shows that $X_i+Y'$ is an extension
of $X_i$ by $Y/X$ and both $X_i$ and $Y/X$ are finite length
modules generated by $\bold z$. From our first considerations, we
know that $X_i+Y'$ is generated by $\bold z$, thus also $Y =
\bigcup_i (X_i+Y')$ is generated by $\bold z.$
        \bigskip
Let $\bold q$ be a class of indecomposable $\Lambda$-modules of
finite length. We want to find a criterion for $g(\bold q)$ to be
the torsion class of a split torsion pair in Mod $\Lambda$. We
denote by $K_0(\Lambda)$ the Grothendieck group of all finite
length $\Lambda$-modules modulo exact sequences. In case
$\delta\:K_0(\Lambda) \to \Bbb Z$ is an additive map and $M$ is a
finite length module, we will write $\delta(M)$ for the value
taken by $\delta$ on the equivalence class of $M$ in
$K_0(\Lambda).$
        \medskip
We say that the class $\bold q$ of indecomposable
$\Lambda$-modules of finite length is {\it closed under
successors} provided given indecomposable $\Lambda$-modules
$M_1,M_2$ of finite length with $\Hom(M_1,M_2) \neq 0$, then
$M_1\in \bold q$ implies $M_2\in\bold q$.

We also consider the following finiteness condition (F): If $N$ is
a $\Lambda$-module with $\Hom(\bold q,N) = 0$ and has a submodule
$U \subseteq N$ of finite length such that $N/U$ is generated by
$\bold q$, then $N$ is of finite length.

Finally, let us say that $\bold q$  is {\it numerically
determined} provided there exists a function $\delta\:K_0(\Lambda)
\to \Bbb Z$ such that an indecomposable $\Lambda$-module $M$ of
finite length belongs to $\bold q$ if and only if $\delta(M) > 0.$
    \medskip
{\bf Proposition 1.} {\it Let $\Lambda$ be an artin algebra. Let
$\bold q$ be a class of indecomposable modules in $\mod \Lambda$
closed under successors.

\item{\rm (a)}
If $\bold q$ is numerically determined, then $\bold q$ satisfies
the condition {\rm(F)}.

\item{\rm (b)} If $\bold q$ satisfies the condition {\rm(F)}, then
$g(\bold q)$ is the torsion class of a split torsion pair in $\Mod
\Lambda$. The corresponding torsionfree class is $r(\bold q).$
\par}
    \medskip
Proof:
 (a) Assume that $\bold q$ is numerically determined with associated function
$\delta$. Let $N$ be a $\Lambda$-module with $\Hom(\bold q,N)=0$,
and let $U$ be a finite length submodule of $N$ such that $N/U$ is
generated by $\bold q$. Then all submodules $N'$ of $N$ of finite
length satisfy $\delta(N') \le 0$. In particular, we have
$\delta(U) \le 0$ and we choose a finite length submodule $U'$ of
$N$ with $U \subseteq U'$ such that $\delta(U')$ is maximal. We
claim that $U' = N$. Otherwise, $U'/U$ is a proper submodule of
$N/U$, and since $N/U$ is generated by $\bold q$, there is $Q\in
\bold q$ and a map $f\: Q \to N/U$ with image not contained in
$U'/U.$ Let $U''/U = U'/U + f(Q) \subseteq N/U.$ In this way, we
have found a submodule $U''$ of $N$ with $U' \subset U''$ and such
that $U''/U'$ is a non-zero epimorphic image of a module in $\bold
q$ and thus a non-zero direct sum of modules in $\bold q$. But the
latter condition means that $\delta(U''/U') > 0$ and therefore
$\delta(U') < \delta(U'')$, a contradiction to the choice of $U'$.
Hence $U'=N$, and consequently $N$ has finite length.

(b) Since $\add\bold q$ is closed both under extensions and under
factor modules, Lemma 4 asserts that $g(\bold q)$ is closed under
extensions.

Denote as before by $r(\bold q)$ the class of all
$\Lambda$-modules $L$ with $\Hom(\bold q,L) = 0.$ We want to show
that any exact sequence $0 \to X \to Y \to Z \to 0$ with $X\in
g(\bold q)$ and $Z\in r(\bold q)$ splits.

First, consider the case when $Z$ is of finite length. We may
suppose that $Z$ is indecomposable and also that $Z$ is not
projective. If the given map $Y \to Z$ is not split epi, we obtain
a commutative diagram where the lower sequence is the almost split
sequence ending in $Z$ $$ \CD
 0 @>>> X      @>>> Y @>>> Z @>>> 0 \cr
 @.    @VVf V      @VVV     @|       \cr
 0 @>>> \tau Z @>>> E @>>> Z @>>> 0 \cr
\endCD
$$ Note that $\tau Z$ does not belong to $\bold q$, since $\bold
q$ is closed under successors and $Z$ is not in $\bold q$. But $X$
is generated by $\bold q$, thus we see that $f$ has to be the zero
map. But this implies that the lower sequence splits, which is
impossible.

In order to take care of the case of $ Z $ having arbitrary
length, we show the following: Given a module $Z$ and a chain of
submodules $U_i$ of $Z$ with union $U = \bigcup_i U_i$, then, if
all $Z/U_i$ belong to $r(\bold q)$, also $Z/U$ belongs to $r(\bold
q)$. For the exact sequences $0\to U_i\to Z\to Z/U_i\to 0$ give
rise to the exact sequence $0\to U\to Z\to\lim\limits_\to Z/U_i\to
0$. Since $\Hom(\bold q, Z/U_i)=0$ for all $i$, we have
$\Hom(\bold q,\lim\limits_\to
Z/U_i)\simeq\lim\limits_\to\Hom(\bold q, Z/U_i)=0$, and hence
$\Hom(\bold q, Z/U)=0$.

Now, consider the case of $Z$ being of arbitrary length. We may
suppose that the map $X \to Y$ is an inclusion map. Let $\Cal U$
be the set of submodules $U$ of $Y$ with $X \cap U = 0$ and
$Y/(X+U) \in r(\bold q)$. Since $0$ belongs to $\Cal U$, this set
is non-empty. Given a chain $(U_i)_i$ of elements of $\Cal U$, the
union $U = \bigcup_i U_i$ belongs to $\Cal U$; namely, it is clear
that $X \cap U = 0$; and it follows from above that $Y/(X+U)$
belongs to $r(\bold q)$, since all $Y/(X+U_i)$ belong to $r(\bold
q)$. As a consequence, we may choose a $U$ maximal in $\Cal U$.
Assume $X+U$ is a proper submodule of $Y$. Let $X+U \subset Y'
\subseteq Y$ with $Y'/(X+U)$ being simple. Let $Y''/Y'$ be the
largest submodule of $Y/Y'$ generated by $\bold q$. Since $g(\bold
q)$ is closed under extensions, it follows that $Y/Y''$ belongs to
$r(\bold q).$ As a submodule of $Y/(X+U)$ the module $Y''/(X+U)$
belongs to $r(\bold q)$. Condition (F) asserts that $Z' =
Y''/(X+U)$ is of finite length. According to the first part of the
proof, we know that $\Ext^1(Z',X) = 0$. Since the embedding $X
\simeq (X+U)/U \subset Y''/U$ has cokernel $Z'$, there exists a
submodule $U'$ of $Y$ containing  $U$ with $(X+U)\cap U' = U$ and
$(X+U)+U' = Y''.$ We see that $X\cap U' = 0$ and that $Y/(X+U') =
Y/Y''$ belongs to $r(\bold q)$, thus $U'$ belongs to $\Cal U$, a
contradiction to the maximality of $U$. Hence $Y = X+U = X\oplus
U$, thus the sequence $0 \to X \to Y \to Y/X \to 0$ splits.
        \medskip
{\bf Remark.} Let us note that these considerations can be
extended as follows: We say that a class $\bold q$   of
indecomposable modules of finite length is {\it numerically almost
determined} provided there exists a function $\delta\:K_0(\Lambda)
\to \Bbb Z$ with the following properties: (i) If $M$ belongs to
$\bold q$, then $\delta(M) \ge 0$, and $\delta(M) > 0$ for all but
a finite number of isomorphism classes of modules $M$ in $\bold
q$; (ii) any indecomposable $\Lambda$-module $M$ in $\mod \Lambda$
with $\delta(M) > 0$ belongs to $\bold q$. We claim: {\it If
$\Lambda$ is an artin algebra and $\bold q$ is a class of
indecomposable modules in $\mod \Lambda$ which is closed under
successors and numerically almost determined, then $\bold q$
satisfies the condition} (F).

Proof:  Assume that $\bold q$ is numerically almost determined
with associated function $\delta$. Let $N$ be a $\Lambda$-module
with $\Hom(\bold q,N)=0$.

We first observe that for any submodule $U'$ (of finite length) of
$N$ there is a bound $b$ with the following property: If $U''$ is
a submodule of $N$ of finite length with $U'\subseteq U''$ such
that $U''/U'$ is generated by $\bold q$ and $\delta(U''/U') = 0$,
then the length of $U''/U'$ is bounded by $b$. Namely, let
$Q_1,\dots,Q_m$ be the indecomposable modules in $\bold q$ (one
from each isomorphism class) with $\delta(Q_i) = 0$ for $1\le i
\le m$, and assume that these modules $Q_i$ are of length at most
$d$. Let $\dim_k\Ext^1(Q_i,U') \le c$ for all $i$. Assume $U'
\subseteq U'' \subseteq N$ is given with $U''$ of finite length,
$U''/U'$ generated by $\bold q$ and $\delta(U''/U') = 0$. If we
write $U''/U'$ as a direct sum of indecomposables, all these
direct summands $X$ belong to $\bold q$ (since $\bold q$ is closed
under successors), thus $\delta(X) \ge 0$. But these numbers add
up to zero, thus we have $\delta(X) = 0$. This shows that $U''/U'
\simeq \bigoplus_i Q_i^{t_i},$ for some natural numbers $t_i$. Now
if $t_i > c$ for some $i$, then $\dim_k\Ext^1(Q_i,U') \le c$
implies that $U''$ has a submodule isomorphic to $Q_i$, in
contrast to the fact that $\Hom(Q_i,N) = 0.$ Altogether we see
that the length of $U''/U'$ is bounded by $b = cd.$

Now, let $U$ be a finite length submodule of $N$ such that $N/U$
is generated by $\bold q$ and choose (as in the proof of
Proposition 1) a finite length submodule $U'$ of $N$ with $U
\subseteq U'$ such that $\delta(U')$ is maximal. Consider chains
$U' = U_0 \subseteq U_1 \subseteq U_2 \subseteq \dots$ of finite
length modules such that the factors $U_i/U_{i-1}$ are generated
by $\bold q$ and satisfy $\delta(U_i/U_{i-1}) = 0$ for all $i$. We
claim that such a sequence stabilizes. For $U_i/U_0$ is generated
by $\bold q$ and $\delta(U_i/U_0) = 0$, thus, as we have seen,
$U_i/U_0$ is of bounded length, with a bound only depending on
$U_0$. It has the following consequence: Replacing $U'$ if
necessary by a larger submodule, we may assume in addition that
any finite length submodule $U''$ of $N$ with $U' \subseteq U''$
which is generated by $\bold q$ satisfies $\delta(U''/U') > 0.$
This then allows to complete the proof as above.
    \bigskip\bigskip
%======================================================================
{\bf 2\. (Concealed) canonical algebras and separating tubular
families.}
        \medskip
In this section we recall some background material on canonical
and concealed canonical algebras.
        \medskip
Given a class $\bold x$ of indecomposable modules of finite
length, we say that an indecomposable module $M$ of finite length
is a {\it proper predecessor of $\bold x$} provided it does not
belong to $\bold x$, but there is a sequence of indecomposables $M
= M_0, M_1,\dots, M_n$ with $\Hom(M_{i-1},M_i) \neq 0$ for all $1
\le i \le n$ such that $M_n$ belongs to $\bold x$. Similarly, $M$
is said to be {\it a proper successor of $\bold x$} provided it
does not belong to $\bold x$, but there is a sequence of
indecomposables $M_0, M_1,\dots, M_n = M$ with $\Hom(M_{i-1},M_i)
\neq 0$ for all $1 \le i \le n$ such that $M_0$ belongs to $\bold
x$.
        \medskip
{\bf Separating tubular families} (see [R2,R4,LP,RS]). Let
$\Lambda$ be an artin algebra, and let $\bold t$ be a sincere
stable separating tubular family. Recall that this means the
following: a {\it tubular family} consists of all the
indecomposables belonging to a set of tubes in the
Auslander-Reiten quiver of $\Lambda$ (in particular, all the
modules in $\bold t$ are of finite length). Such a tubular family
is said to be {\it stable} provided all the tubes are stable, thus
provided it does not contain any indecomposable module which is
projective or injective. A family of modules is said to be {\it
sincere} provided every simple $\Lambda$-module occurs as the
composition factor of at least one of the given modules. Finally,
let us say that the tubular family $\bold t$ is {\it separating}
provided it is standard, there are no indecomposable modules $M$
of finite length which are both proper predecessors of $\bold t$
and proper successors of $\bold t$, and any map from a proper
predecessor of $\bold t$ to a proper successor of $\bold t$
factors through any of the tubes in $\bold t$.

Now let $\bold t$
be a separating tubular family. We denote by $\bold p$ the class
of indecomposables of finite length which are proper predecessors
of $\bold t$, and by $\bold q$ the class of indecomposables of
finite length which are proper successors of $\bold t$. Then any
indecomposable module of finite length belongs either to $\bold
p,$ $\bold t$ or $\bold q$, 
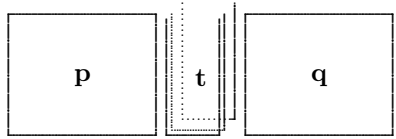
$$ \hbox{\beginpicture
\setcoordinatesystem units <.7cm,.7cm> \put{} at 0 0 \put{} at 6
2.5 \plot 0 0   2.8 0  2.8 2.3    0 2.3  0 0 / \plot 4.5 0 7.3 0
7.3 2.3  4.5 2.3  4.5 0 / \plot 3 2.2  3 0  4 0  4 2.2 / \plot 4.1
0.1  4.1 2.3 / \plot 4.2 0.3  4.3 0.3  4.3 2.5 / \setdots <1pt>
\plot 3.1 2.3 3.1 0.1  4.1 0.1 / \setdots <2pt> \plot 3.3 2.5 3.3
0.3  4.2 0.3 / \put{$\bold p$} at 1.4 1.1 \put{$\bold t$} at 3.65
1.1 \put{$\bold q$} at 5.9 1.1
\endpicture}
$$ and one says that $\bold t$ {\it separates $\bold p$ from
$\bold q$.} Note that there are no maps ''backwards'': $$
 \Hom(\bold t,\bold p) = \Hom(\bold q,\bold p) =\Hom(\bold q,\bold t) = 0
$$ and any map from a module in $\bold p$ to a module in $\bold q$
can be factored through a module in $\bold t$ (even through one
lying in a prescribed tube inside $\bold t$). In case $\bold t$ is
in addition sincere and stable, then all the indecomposable
projective modules belong to $\bold p$, the indecomposable
injective modules belong to $\bold q$. As a consequence, in this
case the modules which belong to $\bold p$ or $\bold t$ have
projective dimension at most 1, those which belong to $\bold t$
or $\bold q$ have injective dimension at most 1. Also note that a
stable separating tubular family $\bold t$ always is an exact
abelian subcategory of $\mod\Lambda$ (and all the indecomposables
in $\bold t$ are serial when considered as objects in this
subcategory).

The algebras $\Lambda$ with a sincere stable separating tubular
family are the {\it concealed canonical algebras}. They have been
studied in [LM], [LP] and [RS], and we are going to review the
main steps of the construction at the end of this section. The
essential ingredient for many of our considerations are the defect
functions on the Grothendieck groups $K_0(\Lambda).$
        \medskip
{\bf The construction of the canonical algebras} (see [R4]). Let
$k$ be a field. We start with a tame bimodule ${}_FM_G$, thus
$F,G$ are division rings having $k$ as central subfield and being
finite-dimensional over $k$, and $M$ is an $F$-$G$-bimodule with
$\dim {}_FM \cdot \dim M_G = 4$ and such that $k$ operates
centrally on $M$. This means that $\Lambda_0= \left[\smallmatrix
F&M\cr 0 &G\endsmallmatrix\right]$ is a finite-dimensional tame
hereditary $k$-algebra with precisely two simple modules, and, up
to Morita equivalence, all finite-dimensional tame hereditary
$k$-algebras with precisely two simple modules are obtained in
this way. A non-zero $\Lambda_0$-module $N$ is called {\it simple
regular,} provided $\tau N \simeq N$ and $\End(N)$ is a division
ring. It is well-known that there are many simple regular
$\Lambda_0$-modules, the number of isomorphism classes is
$\max(\aleph_0,|k|)$.

If $R$ is any ring, $N$ any $R$-module with endomorphism ring
$D^{\text{op}}$ and $n\ge 1$ a natural number, we denote by
$R[N,n]$ the $n$-point extension of $R$ by $N$, it is the matrix
ring $$ R[N,n] \ = \
 \bmatrix R & N & \cdots & N \cr
          0 & D & \cdots & D \cr
          \vdots & \ddots & \ddots & \vdots \cr
          0 & \cdots & 0 & D \endbmatrix.
$$ Since any $R$-module may be considered (in a natural way) as an
$R[N,n]$-module, we may iterate this procedure: given a finite
sequence $N_1,\dots,N_t$ of $R$-modules and natural numbers
$n_1,\dots,n_t$, we may form $R[N_1,n_1]\cdots [N_t,n_t]$.

Let us return to $\Lambda_0$. Choose $t$ pairwise non-isomorphic
simple regular $\Lambda_0$-modules $N_1,\dots,N_t$ (with
endomorphism rings $D_i^{\text{op}}$) and natural numbers
$n_1,\dots,n_t$ and consider $\Lambda' = \Lambda_0[N_1,n_1]\cdots
[N_t,n_t]$, a so called {\it squid algebra}; its quiver (or better
species) is of the form as shown to the left: $$
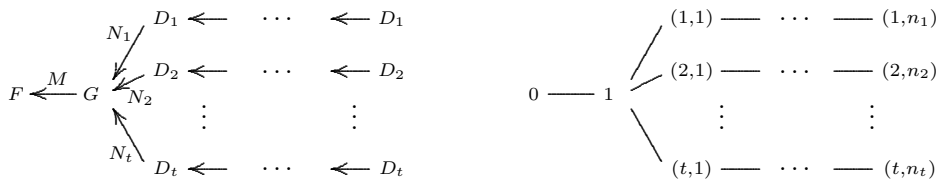
\hbox{\beginpicture \setcoordinatesystem units <1cm,1cm>
\put{\beginpicture \setcoordinatesystem units <1cm,1cm> \put{} at
0 -1 \put{} at  5 1 \put{$\ssize F$} at 0 0 \put{$\ssize G$} at 1
0 \put{$\ssize D_1$} at 2 1 \put{$\ssize D_1$} at 5 1 \put{$\ssize
D_2$} at 2 0.3 \put{$\ssize D_2$} at 5 0.3 \put{$\ssize D_t$} at 2
-1 \put{$\ssize D_t$} at 5 -1 \put{$\ssize M$} at 0.55 0.2
\put{$\ssize N_1$} at 1.37 0.8 \put{$\ssize N_2$} at 1.65 -0.05
\put{$\ssize N_t$} at 1.4 -.8

\arrow <2mm> [0.25,0.75] from .8 0     to  .2 0 \arrow <2mm>
[0.25,0.75] from 1.7 0.9  to 1.3 0.2 \arrow <2mm> [0.25,0.75] from
1.7 -.9  to 1.3 -.2 \arrow <2mm> [0.25,0.75] from 1.7 .25  to 1.3
0.05 \arrow <2mm> [0.25,0.75] from 2.8 1    to 2.3 1 \arrow <2mm>
[0.25,0.75] from 4.7 1    to 4.2 1 \arrow <2mm> [0.25,0.75] from
2.8 .3   to 2.3 .3 \arrow <2mm> [0.25,0.75] from 4.7 .3   to 4.2
.3 \arrow <2mm> [0.25,0.75] from 2.8 -1   to 2.3 -1 \arrow <2mm>
[0.25,0.75] from 4.7 -1   to 4.2 -1

\put{\dots} at 3.5 1 \put{\dots} at 3.5 .3 \put{\dots} at 3.5 -1
\put{$\vdots$} at 2.5 -0.2 \put{$\vdots$} at 4.5 -0.2
\endpicture} at 0 0
\put{\beginpicture \setcoordinatesystem units <1cm,1cm> \put{} at
0 -1 \put{} at  5 1 \put{$\ssize 0$} at 0 0 \put{$\ssize 1$} at 1
0 \put{$\ssize (1,1)$} at 2.1 1 \put{$\ssize (1,n_1)$} at 5 1
\put{$\ssize (2,1)$} at 2.1 0.3 \put{$\ssize (2,n_2)$} at 5 0.3
\put{$\ssize (t,1)$} at 2.1 -1 \put{$\ssize (t,n_t)$} at 5 -1
\plot 0.2 0 0.8 0 / \plot 1.3 0.2 1.7 0.9 / \plot 1.3 -0.2 1.7
-0.9 / \plot 1.3 0.05 1.7 0.25 / \plot 2.5 1 3 1 / \plot 4 1 4.5 1
/

\plot 2.5 .3 3 .3 / \plot 4 .3 4.5 .3 /

\plot 2.5 -1 3 -1 / \plot 4 -1 4.5 -1 / \put{\dots} at 3.5 1
\put{\dots} at 3.5 .3 \put{\dots} at 3.5 -1 \put{$\vdots$} at 2.5
-0.2 \put{$\vdots$} at 4.5 -0.2
\endpicture} at 7 0
\endpicture}
$$ As shown to the right, we label the two vertices of $\Lambda_0$
by $0$ and $1$, and the extension vertices of the $i$-th branch by
$(i,1),\dots,(i,n_i)$, always from left to right. It is not
difficult to see that $I(0)\oplus I(1)\oplus \bigoplus_{i,j}
\tau^{j} I(i,j)$ is a cotilting module; its endomorphism ring is
denoted by $\Lambda$ and the algebras $\Lambda$ obtained in this
way are the {\it canonical} algebras.

Given the canonical algebra $\Lambda$, let us write down the
(canonical) defect function $\delta\:K_0(\Lambda) \to \Bbb Z$:
Note that $\Lambda$ has a unique simple projective module $S$ and
a unique simple injective module $S'$. The defect $\delta(M)$ of a
$\Lambda$-module $M$ is calculated in terms of the Jordan-H\"older
multiplicities $[M:S]$ and $[M:S']$, as follows: $$
 \delta(M) = \left\{\matrix\quad \phantom{2}[M:S']-\phantom{2}[M:S] &&\dim_k S' = \dim_k S \cr
                           \quad2[M:S']-\phantom{2}[M:S] &\text{ in case }&\dim_k S' > \dim_k S \cr
                           \quad\phantom{2}[M:S']-2[M:S] &&\dim_k S' < \dim_k S \endmatrix
                           \right.
$$ If we denote by $\bold t$ the class of all indecomposable
$\Lambda$-modules $M$ with $\delta(M) = 0$, then $\bold t$ is a
stable separating tubular family, separating the class $\bold p$
of all indecomposable modules $M$ with $\delta(M) <0$ from the
class $\bold q$ of all indecomposable modules $M$ with $\delta(M)
> 0.$ We call this triple $(\bold p,\bold t,\bold q)$ the {\it
canonical trisection} of $\mod \Lambda$.
        \medskip
{\bf The construction of the concealed canonical algebras.} Let
$\Lambda$ be a canonical algebra, with canonical trisection
$(\bold p,\bold t,\bold q)$. Let $T$ be a tilting module which
belongs to $\add \bold p$ (since all the modules in $\bold p$ have
projective dimension at most 1, to be a tilting module means in
addition that $\Ext^1(T,T) = 0$ and that there is an exact
sequence $0 \to {}_\Lambda\Lambda \to T' \to T'' \to 0$ with
$T',T'' \in \add T$). Then, by definition, $\Lambda' =
\End(T)^{\text{op}}$ is a {\it concealed canonical} algebra. We
note that the tilting functor $F = \Hom(T,-)$ sends $\bold t$ to a
sincere stable separating tubular family $\bold t'$ in $\mod
\Lambda'$, and as it has been shown in [LM] and [LP] all sincere
stable separating tubular families are obtained in this way.
         \bigskip\bigskip
%===================================================
{\bf 3\. Two extremal torsion pairs in $\Mod\Lambda$}
    \medskip
In the next sections, let $\Lambda$ be a canonical algebra and
$(\bold p,\bold t,\bold q)$ its canonical trisection. In this
section we introduce two torsion pairs in $\Mod \Lambda$ which
will turn out to be split, both having the property that the
indecomposable torsion modules of  finite length are just the
modules in $\bold q$. We also introduce and investigate the class
of Pr\"{u}fer modules.
    \medskip
{\bf The torsion pair $(\Cal C,\Cal Q)$.} As we have mentioned,
the category $\bold q$ is closed under successors. Since it is
also numerically determined, we are able to apply Proposition 1.
Thus, if we denote $\Cal C = r(\bold q)$ and $\Cal Q = g(\bold
q)$, then $(\Cal C, \Cal Q )$  is a split torsion pair in $\Mod
\Lambda$.
    \medskip
%==============================================
{\bf The torsion pair $(\Cal R,\Cal D).$} Let $\Cal D = l(\bold
t)$; note that $\Cal D$ can also be described as $\Cal D = \{M\mid
\Ext^1(\bold t,M) = 0\}.$ Namely, the objects $T$ in $\bold t$
have projective dimension 1, thus we have $\Ext^1(T,M) \cong
D\Hom(M,\tau T)$ (here, $D = \Hom_k(-,k)$ is the duality with
respect to the base field $k$). Since $\bold t$ consists of stable
tubes, the Auslander-Reiten translation is bijective on the
isomorphism classes in $\bold t$.

Since $\Cal D = l(\bold t)$, it is the torsion class of a torsion
pair in  $ \Mod \Lambda$, namely of $ (\Cal R, \Cal D)$, where
$\Cal R = r(\Cal D) = rl(\bold t)$ is the smallest torsionfree
class containing the class $\bold t$. As we have mentioned above,
we may describe $\Cal R$ also as follows: A module $M$ belongs to
$\Cal R$ if and only if the only submodule $U$ of $M$ with
$\Hom(U,\bold t) = 0$ is $U = 0.$ (We later will see that also
this torsion pair is split.)
        \bigskip
The two torsion pairs $(\Cal C,\Cal Q)$ and $(\Cal R,\Cal D)$ are
related as follows: $$ \Cal R \subseteq \Cal C \t{and}
 \Cal Q  \subseteq \Cal D.
$$ These two assertions are equivalent, thus it is sufficient to
verify one of them. But actually, both follow directly from the
assertion $\Hom(\bold q,\bold t) = 0.$
    \medskip
There is the following straightforward characterization of these
two torsion pairs as the extremal ones when dealing with all the
torsion pairs of $\Mod \Lambda$ with prescribed distribution of
the finite dimensional modules: we deal with the torsion pairs
$(\Cal X,\Cal Y)$ with $\bold t\subset \Cal X$ and $\bold q
\subset \Cal Y$.
    \medskip
{\it The two torsion pairs $(\Cal R,\Cal D)$ and $(\Cal C,\Cal Q
)$ have the property that the finite dimensional indecomposable
torsionfree modules are those in $\bold p$ and $\bold t$, whereas
the finite-dimensional indecomposable torsion modules are those in
$\bold q$. If $(\Cal X,\Cal Y)$ is an arbitrary torsion pair in
 $\Mod\Lambda$ such that the finite dimensional
indecomposable torsionfree modules are those in $\bold p$ and
$\bold t$, or, equivalently, such that the finite dimensional
indecomposable torsion modules are those in $\bold q$, then} $$
 \Cal R \subseteq \Cal X \subseteq \Cal C \t{and}
 \Cal Q  \subseteq \Cal Y \subseteq \Cal D.
$$
    \medskip
%==============================================
{\bf The intersection $\omega$ of $\Cal C$ and $\Cal D$.} Let
$\omega = \Cal C \;\cap\; \Cal D. $ 
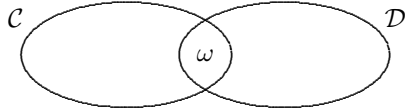
$$ \hbox{\beginpicture
\setcoordinatesystem units <.7cm,.7cm> \put{} at  0 0 \put{} at  7
2 \ellipticalarc axes ratio 2:1 360 degrees  from 0 1 center at 2
1 \ellipticalarc axes ratio 2:1 360 degrees from 3 1 center at 5 1
\put{$\omega$} at 3.5 1 \put{$\Cal C$} at -0.1 1.7 \put{$\Cal D$}
at 7.1 1.7
%\put{$c(\bold t)$} at 1.5 0.7
%\put{$\Cal Q $} at 5.5 0.7
\endpicture}
$$ A complete description of $\omega$ and its relation to $\Cal C$
and $\Cal D$ is one of the aims of the paper.
    \medskip
%=======================================================
{\bf The torsion modules $\Cal T$.} Let $\Cal T =  \Cal C \cap
g(\bold t)$, these are the modules in $\Cal C$ generated by $\bold
t$. We stress that {\it every module in $\Cal T$ is the union of
modules in $\add\bold t$ and that $\Cal T = \lim\limits_\to \bold
t$,} the direct limit closure of $\add\bold t$. Namely, if $M =
\lim\limits_\to M_i$ is the direct limit of a directed system of
modules $M_i$ in $\add\bold t$, then $M$ obviously belongs to
$g(\bold t)$;  since $\Hom(Y,M) = \lim\limits_\to \Hom(Y,M_i)$ for
any finitely generated module $Y$, we also see that $\Hom(Y,M) =
0$ for $Y \in \bold q,$ thus $M$ also belongs to $r(\bold q) =
\Cal C.$ And conversely, if we assume that $M$ is generated by
$\bold t$, then $M$ is the union of all its submodules $M_i$ which
are images of maps from modules in $\add\bold t$ to $M$, and this
is a directed system. As a factor module of a module in $\add\bold
t$, any $M_i$ is a direct sum of modules in $\bold t$ and in
$\bold q$. If we assume in addition that $M$ belongs to $r(\bold
q)$, then all the $M_i$ belong to $\add\bold t.$ This shows that
$M$ is the union of a directed system of submodules which belong
to $\add\bold t.$

Since $\bold t$ is an exact abelian subcategory of $\mod\Lambda$,
it follows that {\it $\Cal T$ is an exact abelian subcategory of
$\Mod\Lambda$.} In particular, $\Cal T$ is closed under kernels,
images and cokernels, and also under direct sums.
        \medskip
%%%%%%%%%%%%%%%%HERE: \omega_0%%%%%%%%%%%%%%%%%%%%%%%%%
%====================================================
{\bf The Pr\"ufer modules.} Of special interest are the so called
{\it Pr\"ufer modules}. They are constructed as follows: The full
subcategory $\add \bold t$ given by the finite direct sums of
modules in $\bold t$ is an abelian length category. Every module
in $\bold t$ has a unique composition series when considered
inside this subcategory $\add \bold t$; its length is called the
{\it regular length}; the modules in $\bold t$ of regular length 1
are just the simple objects of $\add \bold t$ (we call them the
{\it simple objects} of $\bold t$).

The isomorphism classes in $\bold t$ are indexed by pairs
consisting of a natural number $r$ and (the isomorphism class of)
a simple object $S$ in $\bold t$. We will denote the corresponding
$\Lambda$-module by $S[r]$. It is the unique module in $\bold t$
with regular length $r$ and having $S$ as a submodule. For any
simple object $S$ in $\bold t$, there is a sequence of inclusion
maps $$
 S = S[1] \to S[2] \to \cdots \to S[r] \to \cdots,
$$ and we denote by $S[\infty]$ the direct limit of this sequence.
This is the {\it Pr\"ufer module with regular socle $S$}.

We note the following: As an object in $\Cal T$, any Pr\"ufer
module $S[\infty]$ belongs to $\Cal C$. Since $S[\infty]$ is
injective in $\Cal T$, we have in particular $\Ext^1(T,S[\infty])
= 0$ for any object $T$ in $\bold t$, thus $S[\infty]$ also
belongs to $\Cal D$. This shows: {\it The Pr\"ufer modules belong
to $\omega$.} We can strengthen this assertion as follows: {\it
The direct sums of Pr\"ufer modules are just the injective objects
of the abelian category $\Cal T$, and every object in $\Cal T$ has
an injective envelope} (see [R1]). Let us denote $\omega_0 = \Cal
T \cap \Cal D$, then this is the full subcategory of all injective
objects of $\Cal T$. Thus $\omega_0$ is the full subcategory of
all direct sums of Pr\"ufer modules.
%================================================
    \bigskip
Given a module $M$, we denote by $tM$ the maximal submodule of $M$
generated by $\bold t,$ thus $t(M/tM) = 0$ for any module $M$. We
use Lemma 4 in order to see that the class of modules $ g(\bold t)
$ generated by $ \bold t $ is closed under extensions. As a
consequence, the pair $(\Cal F, g(\bold t)) $ where $\Cal F =
r(\bold t)$, is a torsion pair (but this is a torsion pair which
is not split).
        \bigskip
{\bf The analogy.} The notation introduced above should remind the
reader of the analogous situation when dealing with the category
$\Mod R$, where $R$ is a Dedekind ring, say a Dedekind ring with
infinitely many maximal ideals. The torsion pair in $\Mod R$ which
we have in mind is the usual one: the torsion modules are those
$R$-modules $M$ where every element is annihilated by some
non-zero ideal, the torsionfree $R$-modules are those with zero
torsion submodule. In our situation of dealing with a canonical
algebra $\Lambda$, we consider the torsion pair $(\Cal F,g(\bold
t))$; note that it is the subcategory $\Cal C$ of $\Mod \Lambda$
which shows strong similarity to $\Mod R$, thus we reserve the
symbol $\Cal T$ for the intersection of $g(\bold t)$ with $\Cal
C$. The module class $\Cal D$ should be interpreted as the
``divisible'' modules, the module class $\Cal R$ as the
``reduced'' ones.
    \bigskip\bigskip
%======================================================================
{\bf 4\. The $\omega$-coresolution of the modules in $\Cal C$.}
    \medskip
In this section we show that there are $\omega$-coresolutions for
the modules in $\Cal C$, with $\omega=\Cal C\cap\Cal D$ as before.
        \medskip
{\bf Theorem 1.} {\it For every $\Lambda$-module $M$, there exists
a minimal left $\omega$-approximation, $M \to M_\omega$, and its
cokernel belongs to $\omega_0.$ This minimal left
$\omega$-approximation is injective if and only if $M$ belongs to
$\Cal C$.}

{\it If $M$ belongs to $\Cal F$, then $M_\omega$ belongs to $\Cal
F.$ If $M$ belongs to $\Cal T$, then $M_\omega$ belongs to $\Cal
T.$}
    \medskip
Part of the theorem may be reformulated as follows: {\it For any
$M \in \Cal C$, there is an exact sequence $$
 0 @>>> M @>f >> M_\omega @>>> T @>>> 0
$$ with $M_\omega\in \omega$ and $T\in \omega_0$, such that $f$ is
a minimal left $\omega$-approximation.} In this way, one obtains a
characterization of the modules in $\Cal C$ as follows: {\it The
modules in $\Cal C$ are the kernels of epimorphisms in $\omega$.}
        \medskip
Since $(\Cal C,\Cal Q) $ is a split torsion pair, the module $M$
is a direct sum of a module in $\Cal C$ and a module in $\Cal Q $.
For $M \in \Cal Q $, the minimal left $\omega$-approximation
$M_\omega$ has to be zero, since $\Hom(\bold q,\Cal C) = 0.$ For
the proof of Theorem 1, it is sufficient to construct a minimal
left $\omega$-approximation for the modules in $\Cal C$. First, we
construct an exact sequence $$
 0 \to M \to M_\omega \to T \to 0
$$ where $M_\omega$ is in $\omega$ and $T\in \omega_0$.
    \bigskip
For the proof we will need a splitting result which later will be
incorporated into our basic splitting theorem (Theorem 2):
    \medskip
{\bf Lemma 5.} $$
 \Ext^1(\Cal T,\Cal D) = 0.
$$
    \medskip
This is an immediate consequence of the fact that $\Cal T =
\lim\limits_\to \bold t$.
    \medskip
Proof of Theorem 1. As we have mentioned, we can assume that $M$
belongs to $\Cal C$. Take a universal extension $$
 0 @>>> M @>\mu'>> M' @>\pi'>> T' @>>> 0,
$$ with $T'$ a direct sum of simple objects in $\bold t$ and let
$\epsilon$ be its equivalence class in $\Ext^1(T',M)$. The
universality means the following: given a simple object $S$ in
$\bold t$, then first, any element of $\Ext^1(S,M)$ is induced
from $\epsilon$ by a map $S \to T'$, and second, that $\pi'f = 0$
for any map $f\:S \to M'$. Note that the first of these conditions
can be reformulated as saying that $\Ext^1(S,\mu') = 0.$

Take an injective envelope $u\:T' \to T$ in the abelian category
$\Cal T$, thus $T\in \omega_0$. The cokernel of $u$ belongs to
$\Cal T$, thus it has projective dimension at most 1 since it is a
direct limit of finite length modules of projective dimension 1.
It follows that the map $\Ext^1(u,M)$ is surjective, thus there
exists a commutative diagram with exact rows $$ \CD
 0 @>>> M @>\mu'>> M'       @>\pi'>> T' @>>> 0 \cr
 @.    @|         @Vu'VV           @Vu VV      \cr
 0 @>>> M @>\mu>>  M_\omega @>\pi>>  T @>>> 0  \cr
\endCD
$$ Since $M$ and $T$ belong to $\Cal C$ and $\Cal C$ is closed
under extensions, we see that $M_\omega$ belongs to $\Cal C$.

In order to show that $M_\omega$ belongs to $\Cal D$, it is enough
to show that $\Ext^1(S,M_\omega)=0$ for all simple objects in
$\bold t$ (then clearly $\Ext^1(T'',M_\omega) =0$ for any object
$T''$ in $\bold t$). The maps $\mu,\pi$ yield  an  exact sequence
$$
 \Ext^1(S,M) @>{\Ext^1(S,\mu)}>> \Ext^1(S,M_\omega)    @>{\Ext^1(S,\pi)}>> \Ext^1(S,T),
$$ and the last term is zero, since $T$ is injective in $\Cal T$.
Thus the map $\Ext^1(S,\mu)$ is surjective. However, this map
$\Ext^1(S,\mu)$ factors through  $\Ext^1(S,\mu') = 0$. Thus we
conclude that $\Ext^1(S,M_\omega) = 0.$

It remains to be seen that the map $\mu$ is a minimal left
$\omega$-approximation. According to Lemma 5, we have
$\Ext^1(T,\omega) = 0$, and thus $\mu$ is a left
$\omega$-approximation.

In order to show that $\mu$ is left minimal, we first show that
for a direct sum decomposition $M_\omega = N\oplus N'$ with
$\mu(M) \subseteq N$ we must have $N' = 0.$ Thus, consider such a
direct sum decomposition $M_\omega = N\oplus N'$ with $\mu(M)
\subseteq N$. The cokernel $T$ of $\mu$ is isomorphic to
$N/\mu(M)\oplus N'$. Assume $N'$ is non-zero. Since $T$ and
therefore $N'$ is isomorphic to a direct sum of Pr\"ufer modules,
there is a monomorphism $f\:S \to N'$ with $S$ a simple object of
$\bold t.$ The image of $f$ has to lie in the image of $u'$, thus
there is $f'\:S \to M'$ with $f = u'f'$. By construction of the
universal extension $\epsilon$, the composition $\pi'f'$ is zero,
thus $f' = \mu'f''$ for some $f''\:S \to M$. But this implies that
the image of $f = u'f' = u'\mu'f'' = \mu f''$ lies in $\mu(M)
\subseteq N$ and not in $N'$. This contradiction shows that $N'=
0$.

Now consider a map $g\:M_\omega \to M_\omega$ with $g\mu = \mu.$
We obtain a commutative diagram $$ \CD
 0 @>>> M @>\mu>> M_\omega       @>\pi>> T @>>> 0 \cr
 @.    @|         @Vg VV           @Vg' VV      \cr
 0 @>>> M @>\mu>>  M_\omega @>\pi>>  T @>>> 0  \cr
\endCD
$$ Note that $\pi$ induces an isomorphism between the kernel of
$g$ and the kernel of $g'$. In order to show that $g$ is
injective, let us assume to the contrary that $\Ker g \simeq \Ker
g'$ is non-zero. Note that the kernel $\Ker g'$ of $g'$ belongs to
$\Cal T$, thus there is a simple object $S$ of $\Cal T$ which is
contained in $\Ker g'$ and hence $S \subseteq T'$. The isomorphism
of kernels $\Ker g \simeq \Ker g'$ shows that this $S$ may be
considered as a submodule of $M'$ with non-zero composition $S \to
M' \to T'$, but this is a contradiction. Thus $g$ is a
monomorphism.

To see that $g$ is an epimorphism, denote by $L$ its cokernel
$\Cok g.$ Also for the cokernels, $\pi$ induces an isomorphism
$\Cok g \to \Cok g'$. Since $g'\:T \to T$ is a split monomorphism,
and $T\in \omega_0$, we see that $L\in \omega_0$, thus belongs to
$\omega$. Using $\Ext^1(\omega,\omega) = 0,$ we see that $g$ is a
split monomorphism. But according to the previous considerations
this implies that $g$ is surjective.

Of course, if $M$ belongs to $\Cal T$, then also $M_\omega$
belongs to $\Cal T$, since $\Cal T$ is closed under extensions.
Thus, finally, consider the case when $\Hom(\bold t,M) = 0.$ In
order to show that $\Hom(\bold t,M_\omega) = 0$, it is sufficient
to show that $\Hom(S,M_\omega) = 0$ for any simple object $S$ in
$\bold t$. Thus, take a nonzero map $f\:S \to M_\omega$. Its
composition with $\pi$ goes to the socle of $T$ in $\Cal T$, thus
$f = u'f'$ for some $f'\:S \to M'$. But by construction $\pi'f' =
0,$ thus $f' = \mu'f''$ for some $f''\:S \to M$. Since we assume
that there are no non-zero maps $S \to M$, it follows that $f =
0.$ Hence $M_\omega$ is in $\Cal F$.
        \medskip
{\bf Lemma 6.} {\it Assume that $M$ belongs to $\Cal C$ and is of
finite length. If we write $M_\omega/M$ as a direct sum of
Pr\"ufer modules, then any Pr\"ufer module occurs with finite
multiplicity.}
        \medskip
Proof: We may assume that $M$ is indecomposable. If $M$ belongs to
$\bold t$, then $M_\omega$ and $M_\omega/M$ are Pr\"ufer modules
themselves.

Thus, we may assume that $M$ belongs to $\bold p$ and therefore to
$\Cal F$. Let $S$ be a simple object of $\bold t$ and $S[\infty]$
the corresponding Pr\"ufer module. Let $n = \dim_k\Ext^1(S,M)$. We
claim that $S[\infty]$ occurs in $M_\omega/M$ with multiplicity at
most $n$. Otherwise, $M_\omega/M$ has a submodule of the form
$S^{n+1}$, say $U/M$, where $U$ is a submodule of $M_\omega$ with
$M \subseteq U$. Since $\dim_k\Ext^1(S,M) = n$, it follows that
$U$ has a submodule isomorphic to $S$. This is impossible, since
$M_\omega$ belongs to $\Cal F$.
    \bigskip\bigskip
%======================================================================
{\bf 5\. The basic splitting result.}
    \medskip
The aim of this section is to prove the basic splitting result
$\Ext^1(\Cal C,\Cal D)=0$, which will also have as a consequence
that the torsion pair determined by $\Cal D$ splits.

Let $\Lambda$ be a canonical algebra with canonical trisection
$(\bold p,\bold t, \bold q)$ of $\mod \Lambda$ and defect function
$\delta$. We say that an indecomposable projective
$\Lambda$-module $P$ is called a {\it peg} (with respect to
$\delta$) provided $\delta(P) = -1$. If $\Lambda$ is a canonical
algebra, then it is easy to see that a peg exists: if the simple
projective module is not a peg, then the sincere indecomposable
projective module turns out to be a peg. More precisely, let $P$
be the simple projective module and $P'$ the sincere
indecomposable projective module. Then $P$ is a peg if and only if
$\dim\End P \le \dim\End P'$, and $P'$ is a peg if and only if
$\dim\End P \ge \dim\End P'$.
    \medskip
{\bf Lemma 7.} {\it Let $M$ belong to $\Cal \omega$ and let $P$ be
a peg. Then $M$ has a submodule $U$ which is a direct sum of
copies of $P$ such that $M/U$ belongs to $\omega_0.$}
    \medskip
Proof: Let $\Cal U$ be the set of submodules $U'$ of $M$ which are
direct sums of copies of $P$ such that $M/U'$ is in $\Cal C$. We
consider this set as being partially ordered with respect to split
embeddings. Given a chain inside $\Cal U,$ it is not difficult to
see that the union is again in $\Cal U$. Thus we can apply Zorn's
lemma in order to obtain a maximal member $U$ of $\Cal U.$ We show
that $M/U$ belongs to $\Cal T.$ Let $U \subseteq V \subseteq M$
such that $V/U = t(M/U)$. Thus we have to show that $V = M$. The
module $M/V$ belongs to $\Cal F$, and also to $\Cal D$, thus to
$\omega$. The non-zero modules in $\omega$ are sincere. Thus there
is a non-zero homomorphism $f\:P \to M/V$. Note that the kernel
$K$ of $f$ has to be zero, since otherwise $\delta(K) \le
\delta(P)$ and therefore $\delta(P/K) = \delta(P)-\delta(K) \ge
0$. However since $M/V$ is a module in $\Cal F$,
 every non-zero submodule of $M/V$ of finite length
has negative defect. Since $P$ is projective, we can lift the
homomorphism $f\:P \to M/V$ to a homomorphism $f'\:P \to M$ such
that $f = pf',$ where $p\:M \to M/V$ is the canonical map. The
image $P'$ of $f'$ is a submodule of $M$ isomorphic to $P$ and
$P'\cap V = 0$. In particular, we also have $P'\cap U = 0.$ Let
$U' = P'\oplus U$. Then this is a submodule of $M$ which is a
direct sum of copies of $P$. In order to see that the factor
module $M/U'$ belongs to $\Cal C,$ one observes that $M/U'$ is an
extension of the modules $V/U$ and $M/(P'+V).$ The module $V/U$ is
a submodule of $M/U$, thus it belongs to $\Cal C$. If $M/(P'+V)$
would contain a module from $\bold q$ as a submodule, then its
inverse image under the projection $M/V \to M/(P'+V)$ would have
non-negative defect (being an extension of
 the module $(P'+V)/V \simeq P$ of defect $-1$ by a module of positive defect).
But this is impossible, since $M/V$ belongs to $\Cal F$.
Altogether, we see that $U'$ belongs to $\Cal U.$ Since $U$ is a
direct summand of $U'$, we obtain a contradiction to the
maximality of $U$. This shows that $M/U$ belongs to $\Cal T$. As a
factor module of $M\in \Cal D$, the module $M/U$ also belongs to
$\Cal D$, thus to $\omega_0 = \Cal D \cap \Cal T.$
    \bigskip
%=============================================
{\bf Theorem 2 (Basic splitting result):} $$
  \Ext^1(\Cal C,\Cal D) = 0.
$$
    \medskip
Proof. Let $C\in \Cal C$, $D\in\Cal D$. Since $(\Cal C,\Cal Q )$
is a split torsion pair, we can write $D = D'\oplus D''$ with
$D'\in \Cal C$ and $D''\in \Cal Q $. Also, since $(\Cal C,\Cal Q
)$ is a split torsion pair, we have $\Ext^1(\Cal C,\Cal Q ) = 0$.
Thus it is sufficient to show that $\Ext^1(\Cal C,D') = 0$. Note
that $D'$ as a direct summand of $D$ belongs to $\Cal D$, thus to
$\omega$. This shows that we have to show $\Ext^1(\Cal C,\omega) =
0.$

First, we show $\Ext^1(\omega,\omega) = 0.$ Start with a module
$M\in \omega$. According to Lemma 7, there is a submodule $U$
which is a direct sum $\bigoplus_I P$ of copies of a peg $P$ such
that $M/U$ belongs to $\Cal T.$ Let $N$ be a second module in
$\omega$. On one hand, we have $$
 \Ext^1(U,N) = \Ext^1(\bigoplus\nolimits_IP,N) \cong
\prod\nolimits_I \Ext^1(P,N) = 0. $$ On the other hand, Lemma 5
asserts that $\Ext^1(M/U,N) = 0.$ Altogether we conclude that
$\Ext^1(M,N) = 0.$

Now take an arbitrary module $M$ in $\Cal C$, and consider the
minimal left $\omega$-approximation given by Theorem 1: $$
 0 \to M \to M_\omega \to T \to 0
$$ with $M_\omega$ in $\omega$ and $T\in \Cal T$. Applying the
long exact sequence with respect to $\Hom(-,N)$, where $N\in
\omega$, we get the exact sequence $$
 \Ext^1(M_\omega,N) \to \Ext^1(M,N) \to \Ext^2(T,N).
$$ Since the projective dimension of $T$ is at most one, the last
term vanishes. Since $M_\omega$ and $N$ both belong to $\omega$,
also the first term is zero. Thus $\Ext^1(M,N) = 0.$ In this way,
we have shown that $\Ext^1(\Cal C,\omega) = 0,$ as required. This
concludes the proof.
    \bigskip
{\bf Corollary 1.} {\it The torsion pair $(\Cal R,\Cal D)$
splits.}
    \medskip
This follows immediately from the inclusion $\Cal R \subseteq \Cal
C$.
    \bigskip
The category $\Mod \Lambda$ consists of three parts: $$
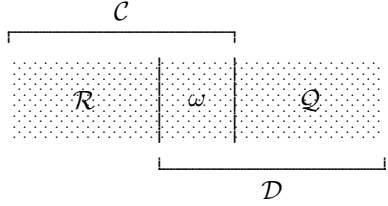
\hbox{\beginpicture \setcoordinatesystem units <1cm,.7cm> \put{}
at  0 0 \put{} at  5 2 \plot 2 0.2  2 1.8 / \plot 3 0.2  3 1.8 /
\plot 0 2.1 0 2.3  3 2.3  3 2.1 / \plot 2 -.1 2 -.3  5 -.3  5 -.1
/ \put{$\omega$} at 2.5 1 \put{$\Cal C$} at 1.5 2.7 \put{$\Cal D$}
at 3.5 -0.7 \put{$\Cal R$} at 1 1 \put{$\Cal Q $} at 4 1
\setshadegrid span <2pt> \hshade 0.2 0 5  1.8  0 5 /
\endpicture}
$$ and this means the following: {\it any $\Lambda$-module $M$ is
a direct sum $M = M_1\oplus M_2\oplus M_3$ with $M_1\in \Cal R,
M_2\in \omega, M_3\in \Cal Q $,  there are no maps ``backwards'':
$$
 \Hom(\omega,\Cal R) = \Hom(\Cal Q ,\Cal R) = \Hom(\Cal Q ,\omega) = 0.
$$ and any map from $\Cal R$ to $\Cal Q$ factors through a module
in $\omega$.} Also the last assertion is an immediate consequence
of previous results: Given a map $h\:M \to N$ with $M \in\Cal R$
and $N \in \Cal Q$, choose a minimal left $\omega$-approximation
$f\:M \to M_\omega$. According to Theorem 1, the map $f$ is
injective and its cokernel $T$ belongs to $\omega_0$. Since
$\Ext^1(T,N) = 0$ (by Theorem 2 or already Lemma 5), we conclude
that $h$ factors through $f$.
        \medskip
{\bf Corollary 2.} {\it We have the following.
\item{\rm(a)} $\pd\,C\leq 1$ for $C$ in $\Cal C$.
\item{\rm(b)} $\id\,D\leq 1$ for $D$ in $\Cal D$.\par}
        \medskip
Proof (a): For $X$ in $\Mod\Lambda$ we have
$\Ext^2(C,X)\simeq\Ext^1(C,\Omega^{-1}X)$. Now $\Omega^{-1}X$ is
generated by injective modules, thus belongs to $\Cal D$. Now
$\Ext^1(\Cal C,\Cal D)=0$ shows that $\Ext^2(C,X) = 0$, hence
$\pd\,C\leq1$.

(b) This follows similarly, using that all the projective modules
belong to $\Cal C$. If $X$ is any $\Lambda$-module, then $\Omega
X$ belongs to $\Cal C$, thus $\Ext^2(X,D) \simeq \Ext^1(\Omega
X,D) = 0.$
    \bigskip\bigskip
%========================================
{\bf 6\. The structure of $\omega.$}
    \medskip
In this section we give the structure of the modules in $\omega$.
    \medskip
{\bf Theorem 3.} {\it Let $P \to P_\omega$ be the minimal left
$\omega$-approximation of a peg $P$ and let $E$ be the
endomorphism ring of $P_\omega$. Then $E$ is a division ring.}
    \medskip
We will denote $P_\omega$ by $G$ and call it the {\it canonical
generic} module (or just the generic module). It will turn out to
be independent of the choice of the peg $P$ and we will
characterize $G$ (up to isomorphism) as the only module in
$\omega$ with endomorphism ring a division ring.
    \medskip
Proof. Let $G = P_\omega$. Note that $G$  contains $P$ as a
submodule and we denote by $p\:G \to T = G/P$ the projection map.
We know that $T$ is a direct sum of Pr\"ufer modules, each
occurring with finite multiplicity, see Lemma 6.

Let $f\:G \to G$ be a nonzero endomorphism. The restriction of $f
$ to $P$ is also non-zero, since otherwise $f$ would yield a
non-zero map $T \to G$. However $G/P$ is a direct sum of
Pr\"{u}fer modules and $G$  belongs to $\Cal F.$ Since
$\delta(P)=-1$, we conclude that the restriction of $f$ to $P$
must be a monomorphism, using that the kernel $\Ker f$ and the
image $f(P)$ are submodules of $G$.

Since $P+f(P)$ is a finitely generated submodule of $G$, there
exists a submodule $P'$ of $G$ with $P+f(P) \subseteq P'$ such
that $P'/P$ belongs to $\bold t.$ Note that we have $\delta(P') =
\delta(P)+\delta(P'/P) = -1$. Since $P'$ is a submodule of $G$, it
has no non-zero submodules of non-negative defect, thus $P'$ has
to be indecomposable.

We claim that $P'/f(P)$ belongs to $\bold t$. The module $P'/f(P)$
has zero defect. If we assume that $P'/f(P)$ does not belong to
$\bold t$, then $P'/f(P)$ has a submodule of positive defect and
its inverse image in $P'$ would yield a non-zero submodule of
non-negative defect, impossible.

Let $T' = G/P'$, with projection map $p'\:G \to T'$. The map $f$
induces maps $f'\:P \to P'$ and $f''\:T \to T'$, thus we deal with
the following commutative diagram $$ \CD
 0 @>>> P @>>>         G @>p>> T  @>>> 0 \cr
 @.     @VV{f'}V     @VVfV         @VV{f''}V   \cr
 0 @>>> P' @>>>         G @>p'>> T'  @>>> 0 \cr
\endCD
$$ and the snake lemma yields an exact sequence $$
 \Ker f' \to \Ker f \to \Ker f'' \to \Cok f' \to \Cok f \to \Cok f'' \to 0.
$$ As we have noted, $\Ker f' = 0.$ Since $f''$ is a map inside
$\Cal T$, its kernel and cokernel both belong to $\Cal T$. Also,
we have shown that the cokernel of $f'$ belongs to $\Cal T$. This
implies that $\Ker f$ and $\Cok f$ belong to $\Cal T$. However $G$
is in $\Cal F$, thus $\Ker f = 0$. This already shows that $f$ is
injective. Also, we see that $\Ker f''$ is a submodule of $\Cok
f'$, thus of finite length.

We claim that $f''$ is surjective. Note that $\add\bold t$ is a
direct sum of serial categories $\add \bold t(\lambda),$ with
$\lambda$ in some index set $\Omega$, and each subcategory $\bold
t(\lambda)$ contains only finitely many isomorphism classes of
simple objects. Of course, $\Cal T$ is a corresponding direct sum
of categories denoted by $\Cal T(\lambda)$ with
$\lambda\in\Omega$. If we decompose $T$ and $T'$ accordingly, we
obtain direct summands $T_\lambda$ of $T$ and $T'_\lambda$ of $T'$
and $f''$ maps $T_\lambda$ into $T_\lambda'$. On the one hand,
$T_\lambda'$ is obtained from $T_\lambda$ by factoring out a
subobject in $\Cal T(\lambda)$ of finite length, thus if we write
$T_\lambda$ and $T'_\lambda$ as direct sums of Pr\"ufer modules,
the numbers of direct summands are equal. On the other hand, $f''$
induces a map $T_\lambda \to T'_\lambda$ with kernel of finite
length. Altogether this implies that $f''$ maps $T_\lambda$ onto
$T'_\lambda.$ Thus $f''$ is surjective.

But $\Cok f'' = 0$ implies that $\Cok f$ is of finite length and
therefore in $\bold t$. Since $\Hom(G,\bold t) = 0$, we see that
$f$ itself is surjective.

In this way, we have shown that any non-zero endomorphism of $G$
is invertible, thus the endomorphism ring of $G$ is a division
ring.
        \medskip
{\bf Corollary 3.} {\it The generic module $G$ is embeddable into
a direct sum of Pr\"ufer modules.}
     \medskip
Proof: Let $P$ be a peg. Choose any Pr\"ufer module, say
$S[\infty]$ (where $S$ is a simple object of $\bold t$). We claim:
{\it The module $P$ is embeddable into $S[\infty]$.} Let $f\:P \to
S[\infty]$ be a non-zero homomorphism. If the kernel $P'$ of $f$
is non-zero, then $\delta(P') \le -1$ and $\delta(P) = -1$. Since
$P/P'$ is an subobject of some module $S[m]$ in $\bold t$, we have
$\delta(P/P') \le 0$, thus $\delta(P/P') = 0.$ Now $P/P'$ is
indecomposable, thus an object in $\bold t$. It follows that
$P/P'$ is a subobject of $S[t]$, where $t$ is the $\tau$-period of
$S$. But $\Hom(P,S[t])$ is finite dimensional, whereas
$\Hom(P,S[\infty])$ is infinite dimensional. This shows that there
are many monomorphisms $P \to S[\infty].$

Starting with a minimal left $\omega$-approximation of $P$ and a
monomorphism $f\:P \to S[\infty]$, we obtain the following
commutative diagram with exact rows:

$$
 \CD
0 @>>> P         @>>> G @>>> T @>>> 0 \\ @.    @VfVV
@Vf'VV  @|      \\ 0 @>>> S[\infty] @>>> N @>>> T @>>> 0 \\
\endCD
$$ Here $T$ and also $S[\infty]$ belong to $\omega_0$. It follows
that the lower sequence splits and that $N = T\oplus S[\infty]$
belongs to $\omega_0$. With $f$ also $f'$ is injective, thus $G$
embeds into an object of $\omega$.
       \medskip
{\bf Remark.} This statement is quite surprising already in the
case of the Kronecker algebra $\Lambda$ (this is the path algebra
of the quiver with two vertices, say $a$ and $b$ and two arrows
starting at $b$ and ending in $a$) where the module category
$\Mod\Lambda$ shows a strong resemblance to the category of all
abelian groups (or better the category of all $k[T]$-modules),
with the canonical generic $\Lambda$-module $G$ corresponding to
$\Bbb Q$ and the Pr\"ufer modules corresponding to the Pr\"ufer
groups. Of course, in sharp contrast to the embedding of $G$ into
a direct sum of Pr\"ufer modules, there does not exist any
embedding of $\Bbb Q$ into a direct sum of Pr\"ufer groups!
    \bigskip
{\bf Theorem 4.} {\it Any module in $\omega$ is a direct sum of
Pr\"ufer modules and of copies of the generic module.}
    \medskip
Proof: Let $M$ be a module in $\omega$. According to Lemma 7,
there is a submodule $U$ which is a direct sum of copies of $P$
such that $M/U$ is a direct sum of Pr\"ufer modules. Let $U \to
U_\omega$ be the minimal left $\omega$-approximation of $U$. We
obtain a commutative diagram as follows: $$ \CD
 0 @>>> U @>>> U_\omega @>>> T  @>>> 0 \cr
 @.     @|     @VVfV         @VVgV   \cr
 0 @>>> U @>>> M        @>>> M/U @>>> 0 \cr
\endCD
$$ The snake lemma yields an isomorphism of the kernel of $f$ and
the kernel of $g$, as well as an isomorphism of the cokernel of
$f$ and the cokernel of $g$. Since $g\:T \to M/U$ is a map in the
abelian subcategory $\Cal T$, the kernel and the cokernel of $g$
both belong to $\Cal T$.

Since $U$ belongs to $\Cal F$, also $U_\omega$ belongs to $\Cal F$
by Theorem 1. Thus the only subobject of $U_\omega$ which belongs
to $\Cal T$ is the zero module. This shows that both $f$ and $g$
are monomorphisms.

Any factor module of a module in $\omega$ belongs to $\Cal D$.
Thus the cokernel $N$ of $g$ belongs to $\Cal T \cap \Cal D
\subset \omega$, and is again a direct sum of Pr\"ufer modules. It
follows that $\Ext^1(N,U_\omega) = 0$, and thus $f$ splits. This
shows that $M$ is the direct sum of $U_\omega$ and $N$. Since $U$
is a direct sum of copies of $P$, we see that $U_\omega$ is a
direct sum of copies of the generic module $G = P_\omega$. Thus
$M$ is a direct sum of Pr\"ufer modules and of copies of $G$.
        \medskip
Note that all the indecomposable modules in $\omega$ have local
endomorphism rings: the endomorphism ring $E$ of the canonical
generic module is a division ring (Theorem 3), the endomorphism
ring of a Pr\"ufer module is a (not necessarily commutative)
discrete valuation ring. As a consequence, the Theorem of
Krull-Remak-Schmidt-Azumaya can be used: the direct sum
decompositions provided in Theorem 4 are unique up to
isomorphisms.
        \medskip
{\bf Corollary 4.} {\it The modules in $\Cal C$ are precisely the
modules cogenerated by $\Cal T$ and also precisely those modules
which can be embedded into a module in $\omega_0$.}
     \medskip
Proof: If a module can be embedded into a module in $\omega$, then
it is cogenerated by Pr\"ufer modules, thus cogenerated by $\Cal
T$. The modules in $\Cal T$ belong to $\Cal C$ and $\Cal C$ is the
torsionfree class of a torsion pair, thus any module cogenerated
by modules in $\Cal C$ belongs to $\Cal C$. It remains to be shown
that any module $M$ in $\Cal C$ can be embedded into a module from
$\omega_0$. Now according to Theorem 1, the module $M$ embeds into
$M_\omega\in \omega$, and according to Theorem 4, we know that
$M_\omega$ is a direct sum of Pr\"ufer modules and of copies of
$G$. We have seen in Corollary 3 that $G$ itself can be embedded
into a direct sum of Pr\"ufer modules, thus $M$ can be embedded
into a direct sum of Pr\"ufer modules.
    \bigskip
{\bf Corollary 5.} {\it If $M$ belongs to $\Cal F$, then
$M_\omega$ is a direct sum of copies of $G$. If $M$ belongs to
$\Cal F$ and has finite length, then $M_\omega$ is a finite direct
sum of copies of $G$.}
    \medskip
Only the last assertion has to be shown. However, given a left
minimal map $X \to Y$ where $X$ is finitely generated and $Y =
\bigoplus_{i\in I} Y_i$ is a direct sum of non-zero modules $Y_i$,
one immediately sees that the index set $I$ has to be finite.
    \bigskip
{\bf Corollary 6.} {\it The module $G$ has finite length as an
$E^{\text{op}}$-module, where $E$ is its endomorphism ring.}
    \medskip
Proof: We may identify $G$ as a vector space with the vector space
$\Hom({}_\Lambda\Lambda,G)$, and this is an identification of
$E^{\text{op}}$-modules. Let us denote the minimal left
$\omega$-approximation of ${}_\Lambda\Lambda$ by $N$. Then this is
a finite direct sum of copies of $G$, say $N \simeq G^n$ for some
$n$. The approximation map ${}_\Lambda\Lambda \to N$ yields an
isomorphism $$
  \Hom({}_\Lambda\Lambda,G) \simeq \Hom(N,G) \simeq \Hom(G^n,G) \simeq E^n,
$$ and all these isomorphisms are isomorphisms of
$E^{\text{op}}$-modules.
    \bigskip\bigskip
%========================================
{\bf 7\. The $\omega$-resolution of the modules in $\Cal D$.}
    \medskip
Using the previous results we can now obtain $\omega$-resolutions
for the modules in $\Cal D$.

\medskip
{\bf Theorem 5.} {\it For every $\Lambda$-module $M$, there exists
a minimal right $\omega$-approximation $M^\omega\to M$.  Its
kernel is a direct sum of copies of the generic module. This
minimal right $\omega$-approximation is surjective if and only if
$M$ belongs to $\Cal D$. If $M$ belongs to $\Cal Q$, then
$M^\omega$ belongs to $\omega_0$.}
    \medskip
Again, we may reformulate the essential part of the theorem: {\it
For any $M \in \Cal Q$, there is an exact sequence $$
 0 @>>> V @>>> M^\omega @>g>> M @>>> 0
$$ with $M^\omega\in \omega_0$ and $V$ a direct sum of copies of
the generic module, such that $g$ is a minimal right
$\omega$-approximation.} And we obtain a characterization of the
modules in $\Cal D$ as follows: {\it The modules in $\Cal D$ are
the cokernels of monomorphisms in $\omega$.}
        \medskip
Proof. If $M$ belongs to $\Cal R$, then $\Hom(\omega,M) = 0$, thus
$0 \to M$ is a minimal right $\omega$-approximation. If $M$
belongs to $\omega$, then the identity map $M \to M$ is a minimal
right $\omega$-approximation.

Thus we may restrict to the case where $M$ belongs to $\Cal Q .$
We claim that in this case $M$ is generated by a direct sum of
Pr\"ufer modules. It is sufficient to show this for a module
$Q\in\bold q$. Since the projective cover of $ Q$ belongs to
$\bold p$ and any map from $\bold p$ to $\bold q$ factors through
$\bold t$, we only have to show that for $S$  a simple object in
$\bold t$  and any natural number $r$ any map $S[r] \to Q$ can be
extended to $S[r+1]$. However, this follows directly from the fact
that $S[r+1]/S[r]$ belongs to $\bold t$ and $\Ext^1(\bold t,\bold
q) = 0.$

Thus there exists an exact sequence $$
 0 @>>> K @>f >> N @>g>> M @>>> 0
$$ with $N\in \omega_0$. Let us show that there exists such a
sequence where $K$ belongs to $\Cal F$. Without loss of
generality, we can assume that the map $f\:K \to N$ is an
inclusion map. Let $tK$ be the maximal submodule of $K$ generated
by $\Cal T$. Since it is the image of a map from a module in $\Cal
T$ to $N\in \Cal T$, it follows that $tK$ belongs to $\Cal T$. If
we factor out $tK$ from $K$ as well as from $N$, we obtain an
exact sequence $$
 0 \to K/tK \to N/tK \to M \to 0,
$$ where $K/tK$ belongs to $\Cal F$. As a factor module of $N$
inside $\Cal T$, the module $N/tK$ belongs to $\omega_0$. So we
can assume that $K$ is in $\Cal F$.

Next, we claim that we even can assume that $K$ is a direct sum of
copies of the generic module. Let $h\:K \to K_\omega$ be the
minimal left $\omega$-approximation of $K$ and form the induced
exact sequence with respect to $h$. We obtain a commutative
diagram of the form $$ \CD
 0 @>>> K      @>f>> N @>>> M @>g>> 0 \cr
 @.    @VVh V      @VVh'V     @|       \cr
 0 @>>>K_\omega @>f'>> N' @>g'>> M @>>> 0 \cr
\endCD
$$ with exact rows. The cokernels of $h$ and $h'$ coincide. Since
by Theorem 1, the cokernel $N''$ of $h'$ belongs to $\omega_0$,
then $N'$ is an extension of $N$ and $N''$ (indeed: a split
extension), thus it belongs to $\omega_0$.

Thus, consider now an exact sequence $$
 0 @>>> K @>>> N @>g>> M @>>> 0
$$ where $K$ is a direct sum of copies of the generic module and
$N\in \omega_0$. Then the map $N \to M$ is a right
$\omega$-approximation, since $\Ext^1(M,K) = 0$ due to the basic
splitting theorem. In order to see that $g$ is right minimal, let
$\eta\:N \to N$ be an endomorphism with $g\eta = \eta,$ thus we
deal with the following commutative diagram $$ \CD
 0 @>>>   K     @>>>   N    @>g>> M    @>>> 0 \cr
 @.     @VV\eta'V    @VV\eta V     @|         \cr
 0 @>>>   K     @>>>   N    @>g>> M    @>>> 0 \cr
\endCD
$$ with exact rows. Since $K$ belongs to $\Add G$, the same is
true for the kernel and the cokernel of $\eta'$. Since $N$ belongs
to $\Cal T$, the same is true for the kernel and the cokernel of
$\eta.$ Thus $\Ker \eta' \simeq \Ker \eta$ belongs to $\Add G \cap
\Cal T = 0$, and also $\Cok\eta' \simeq \Cok\eta$ belongs to $\Add
G \cap \Cal T = 0.$ This shows that $\eta$ is an isomorphism.
        \medskip
{\bf Corollary 7.} {\it Assume that $M$ belongs to $\Cal Q$ and
has finite length. If $g\:M^\omega \to M$ is a minimal right
$\omega$-approximation with kernel $V$, then $V$ is a finite
direct sum of copies of $G$.}
        \medskip
Proof. It is sufficient to show the following: for every finite
length module $M$, the left $E^{\text{op}}$-module $\Ext^1(M,G)$
has finite length, here $E$ is the endomorphism ring of $G$. Let
$\Omega(M)$ be the kernel of a projective cover of $M$. Then
$\Ext^1(M,G)$ is an epimorphic image of $\Hom(\Omega(M),G)$, thus
we want to show that $\Hom(N,G)$ is of finite length as an
$E^{\text{op}}$-module, for every $\Lambda$-module $N$ of finite
length. Choose a free module $F = {}_\Lambda\Lambda^n$ of finite
length which maps onto $N$. Such a map induces an inclusion of
$\Hom(N,G)$ into $\Hom(F,G)$, and $\Hom(F,G)$ is isomorphic to
$G^n$ as an $E^{\text{op}}$-module. Now use Corollary 6.
    \bigskip
As a direct consequence of Theorem 5 we get the following
description of $\Cal C$ and $\Cal D$ in terms of $\omega.$
        \medskip
%========================================
{\bf Proposition 2.} $$ \align
 \Cal C &= \{M \mid \Ext^1(M,\omega) = 0\}, \cr
 \Cal D &= \{M \mid \Ext^1(\omega,M) = 0\}, \cr
\endalign
$$ Proof: We show the first equality: The inclusion $\subseteq$
follows from the basic splitting result. For the inclusion
$\supseteq$, let $\Ext^1(M,\omega) = 0.$ Write $M = N'\oplus N$
with $N'\in \Cal C$ and $N\in \Cal Q $. The minimal right
$\omega$-approximation yields an exact sequence $0 \to K \to
N^\omega \to N \to 0$ where also $K$ is in $\omega$. Since
$\Ext^1(N,\omega) = 0$, we see that $N$ is a direct summand of
$N^\omega$, thus $N$ is in $\omega \subseteq \Cal C.$

The proof of the second assertion uses the corresponding (dual)
arguments.
    \medskip
{\bf Proposition 3.} {\it The class $\omega$ consists of the
relative injective objects inside $\Cal C$ and of the relative
projective objects inside $\Cal D$.} This means: $$
 \omega = \{M \in \Cal C \mid \Ext^1(\Cal C,M) = 0\}
     = \{M\in \Cal D \mid \Ext^1(M,\Cal D) = 0\}.
$$

Proof: That the modules in $\omega$ are relative injective in
$\Cal C$ and relative projective in $\Cal D$ follows again from
the basic splitting result. Conversely, if $M\in \Cal C$ satisfies
$\Ext^1(\Cal C,M) = 0$, then we have in particular
$\Ext^1(\omega,M) = 0$  and therefore $M\in \Cal D$. But $\Cal C
\cap \Cal D = \omega.$ In a similar way, one obtains the second
equality.
        \medskip
{\bf Remark.} We should stress that there is an important
difference between the $\omega$-coresolutions and the
$\omega$-resolutions. As we know, any module $M$ can be written as
$M = M_0\oplus M_1\oplus M_2$ with $M_0 \in \Cal R$, $M_1\in
\omega$ and $M_2$ in $\Cal Q$. Non-trivial left
$\omega$-approximations do exist for modules in $\Cal R$,
non-trivial right $\omega$-approximations for modules in $\Cal Q$.
Whereas the minimal left $\omega$-approximation $M_\omega$ of a
module $M \in \Cal R$ may be an arbitrary module in $\omega$, the
minimal right $\omega$-approximation of a module $M\in \Cal Q$
will be a direct sum of Pr\"ufer modules alone.

There is another substantial difference: let us compare the
possible minimal left $\omega$-approximations $M_\omega$ and
minimal right $\omega$-approximations $M^\omega$ of finite
dimensional modules $M$. Of course, if $M$ belongs to $\bold t$,
then $M_\omega$ is a Pr\"ufer module and $M^\omega = 0$. Consider
the remaining indecomposable modules $M$ of finite length. If $M$
belongs to $\bold p$, then $M_\omega$ is a {\bf finite} direct sum
of copies of $G$ and $M^\omega = 0$. If $M$ belongs to $\bold q$,
then $M^\omega$ is an {\bf infinite} direct sum of Pr\"ufer
modules and $M_\omega = 0.$ If we take into account the cokernel
of the monomorphism $M \to M_\omega$ for $M\in \bold p$ and the
kernel of the epimorphism $M^\omega \to M$ for $M \in \bold q$,
then this strict dichotomy pertains: the cokernel of the
monomorphism $M \to M_\omega$ will be an {\bf infinite} direct sum
of Pr\"ufer modules, the kernel of the epimorphism $M^\omega \to
M$ will be a {\bf finite} direct sum of copies of $G$. But
actually, looking at maps we encounter some astonishing
parallelity: it turns out that both the $\omega$-coresolutions of
the modules in $\bold p$ as well as the $\omega$-resolutions of
the modules in $\bold q$, are maps $X \to Y$, where $X$ is a
finite direct sum of copies of $G$ and $Y$ is an infinite direct
sum of Pr\"ufer modules. For a module $M$ in $\bold p$, we need an
epimorphism of this kind, and $M$ will be the kernel.  For a
module $M$ in $\bold q$, we need a monomorphism of this kind, and
$M$ will be the cokernel.

Let us consider in detail one special example. Let $P$ be a peg,
thus $P_\omega = G$. Consider the $\omega$-coresolution of $P$ $$
 0 \to P \to G \to Y \to 0,
$$ here $Y = G/P$ is  an infinite direct sum of Pr\"ufer modules.
Such an embedding of $P$ into $G$ will remind anyone of the
embedding of $\Bbb Z$ in $\Bbb Q$, with $\Bbb Q/\Bbb Z$ being an
infinite direct sum of Pr\"ufer groups. But note that there does
not exist any embedding of $\Bbb Q$ into a direct sum of Pr\"ufer
groups. In contrast, let $M$ be an indecomposable $\Lambda$-module
in $\bold q$ with projective dimension $\pd M = 1$ and defect
$\delta(M) = -1$. Then it is easy to see (see the proof of
Corollary 7) that $\Ext^1(M,G)$ is one-dimensional as an
$E^{\text{op}}$-space and therefore the $\omega$-resolution of $M$
is of the form $$
 0 \to G \to  Y'  \to M \to 0
$$ with $Y'$ a direct sum of Pr\"ufer modules. Of course, $Y'$ has
to be an infinite direct sum of Pr\"ufer modules. Thus we obtain
an embedding of $G$ into an infinite direct sum of Pr\"ufer
modules, and the cokernel is indecomposable and of finite length.

Given an abelian category $\Cal A$, it is quite customary to form
the quotient category $\Cal A/\Cal A_0$, where $\Cal A_0$ is the
subcategory  of all modules of finite length. In our case, we look
at the quotient category $\Mod \Lambda/\mod\Lambda$. Note that any
finite dimensional module $M$ in $\bold p$, say with
$\omega$-coresolution $M_\omega \to M_\omega/M$, yields an
isomorphism between $M_\omega$ and $M_\omega/M$ in the quotient
category $\Mod \Lambda/\mod\Lambda$. Similarly, any finite
dimensional module $M$ in $\bold q$, say with $\omega$-resolution
$V \to M^\omega$, yields an isomorphism between $V$ and $M^\omega$
in the quotient category $\Mod \Lambda/\mod\Lambda$. In
particular, we obtain in this way isomorphisms  in the quotient
category $\Mod \Lambda/\mod\Lambda$ between $G$ and infinite
direct sums of Pr\"ufer modules.

         \bigskip\bigskip
%========================
{\bf 8\. Further structure of $ \omega $.}
    \medskip
In this section we investigate the maps inside $\omega $ and the
torsion classes in Mod$\Lambda$ with the property that the
indecomposable torsion modules of finite length are just the
modules in $\bold t$.
        \medskip
{\it There are no non-zero maps from a Pr\"ufer module to $G$,}
since a Pr\"ufer module belongs to $\Cal T$ whereas $G$ belongs to
$\Cal F$. On the other hand, {\it any Pr\"ufer module is generated
by $G$.} Namely, let $S$ be a simple object of $\bold t$ and
$p\:M\to S[\infty]$ a projective cover. Since $M$ belongs to $\Cal
F$, its minimal left $\omega$-approximation $M_\omega$ is in $\Cal
F\cap \omega = \Add G$. If we factor $p$ through $M_\omega,$ we
obtain a surjective map $M_\omega \to S[\infty]$. (Actually, the
construction of $G$ shows that $G$ maps onto $(\tau^tS)[\infty]$
for some $t$, but any $(\tau^tS)[\infty]$ maps onto $S[\infty]$.)
 In this way, we
obtain a further characterization of $\Cal D$.
        \medskip
{\bf Corollary 8.} {\it We have $\Cal D = g(G)$.}
        \medskip
Proof: Since $G$ belongs to $\Cal D$ and $\Cal D$ is closed under
direct sums and factor modules, we see that $g(G) \subseteq \Cal
D$. On the other hand, the minimal right $\omega$-approximation of
any module in $\Cal D$ is surjective, according to Theorem 5. Thus
$\Cal D \subseteq g(\omega)$. But any module in $\omega$ is a
direct sum of Pr\"ufer modules and copies of $G$, thus generated
by $G$.
        \bigskip
Also, note that $\add\bold t$ is a direct sum of infinitely many
serial categories $\add \bold t(\lambda),$ with $\lambda$ in some
index set $\Omega$ of cardinality $\max(\aleph_0,|k|).$ Of course,
$\Cal T$ is a corresponding direct sum of categories denoted $\Cal
T(\lambda)$ with $\lambda\in\Omega$. Let us denote by $\Cal
P(\lambda)$ the full subcategory of all direct sums of copies of
the Pr\"ufer modules belonging to $\Cal T(\lambda)$, for any
$\lambda\in \Omega$. Note that for all $\lambda$, $\Cal
T(\lambda)$ contains only finitely many isomorphism classes of
Pr\"ufer modules, and for all but a finite number of $\lambda$
only one.

\medskip
Our three part visualization of $\Mod\Lambda$ can be refined
accordingly: Recall that the modules in $\omega$ are direct sums
of a module in $\Add G$ (the full subcategory of all direct sums
of copies of $G$) and a module in $\omega_0$, and we divide
$\omega_0$ further into the various full subcategories $\Cal
P(\lambda)$. 
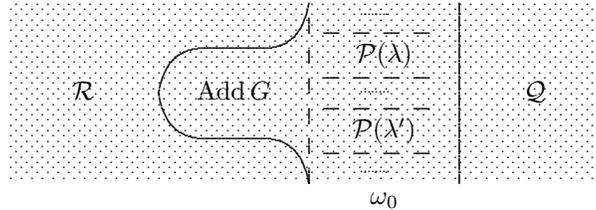
$$ \hbox{\beginpicture \setcoordinatesystem units
<2cm,1cm> \put{} at  0 0 \put{} at  6 2.4
%\plot 2 0  2 2.4 /
\setquadratic \plot 3 0  2.9  0.45 2.7 0.6
        2.5 0.6  2.3 0.6
        2.1 0.75  2   1.2
        2.1 1.65  2.3 1.8
        2.5 1.8  2.7 1.8
        2.9 1.95  3   2.4 /
\setlinear \plot 4 0  4 2.4 / \setdashes <2mm> \plot 3 0  3 2.4 /
\plot 3.1 0.4 3.9 0.4 / \plot 3.1 1   3.9 1 / \plot 3.1 1.4 3.9
1.4 / \plot 3.1 2 3.9 2 / \setdots <.5mm>
%\plot 3.45 0.02  3.45 0.3 /
%\plot 3.45 1.08  3.45 1.36 /
%\plot 3.45 2.1  3.45 2.4 /
\plot 3.35 0.16  3.55 0.16 / \plot 3.35 1.22  3.55 1.22 / \plot
3.35 2.25  3.55 2.25 / \put{$\Add G\strut$} [b] at 2.5 1
\put{$\Cal P(\lambda)$} at 3.5 1.7 \put{$\Cal P(\lambda')$} at 3.5
0.7 \put{$\Cal R\strut$} [b] at 1.5 1 \put{$\Cal Q \strut$} [b] at
4.5 1 \setshadegrid span <2pt> \hshade 0 1 5  2.4  1 5 /
\put{$\omega_0$} at 3.5 -0.22
\endpicture}
$$ Note that the full subcategory $\omega_0$ is separating in the
following sense: First of all, the groups $\Hom(\omega_0,\Cal R)$,
$\Hom(\omega_0,\Add G),$ $\Hom(\Cal Q,\Cal R),$ $\Hom(\Cal Q,\Add
G),$ $\Hom(\Cal Q,\omega_0)$ all are zero, and second, any map
$h\:N \to M$  from a module $N$ in $\Cal R$ or $\Add G$ to a
module $M$ in $\Cal Q$ factors through a module in $\omega_0$
(namely, take a minimal $\omega$-resolution $0 \to V \to M^\omega
\to M \to 0$; since $\Ext^1(N,V) = 0$, the map $h$ factors through
$M^\omega$, but $M^\omega$ belongs to $\omega_0$). The rather
strange shape which we use in order to depict the $\Add G$ part of
$\omega$ should stress that in contrast to $\omega_0$ which is
separating, the subcategory $\Add G$ does not have a corresponding
property.
        \medskip
Using any decomposition of $\Omega$ as the disjoint union of two
subsets $\Omega_1,\Omega_2$, we can write $\add\bold t$ as a
product of two categories, and this yields a corresponding
decomposition of $\Cal T = \Cal T(\Omega_1) \times \Cal
T(\Omega_2)$ as a product of two categories. Any set-theoretical
decomposition $\Omega = \Omega_1\cup \Omega_2$
 therefore gives rise to a split torsion pair
$(\Cal X(\Omega_1),\Cal Y(\Omega_2))$ in $\Mod \Lambda$: the
modules in $\Cal X(\Omega_1)$ are the direct sums $M_1\oplus M_2$,
where $M_1$ is a module in  $\Cal R$ and $M_2$ is the direct sum
of copies of $G$ and of Pr\"ufer modules belonging to $\Cal
T(\Omega_1)$, whereas the modules in $\Cal Y$ are  the direct sums
$M_3\oplus M_4$, where $M_3$ is a direct sum of Pr\"ufer modules
in $\Cal T(\Omega_2)$ and $M_4$ is a module in  $\Cal Q ,$ We have
the following information on torsion pairs.
        \medskip
{\bf Proposition 4.} {\it The only torsion pairs $(\Cal X,\Cal Y)$
with $\bold t \subset \Cal X$ and $\bold q \subset \Cal Y$ are
$(\Cal R,\Cal D)$ and those of the form $(\Cal X(\Omega_1),\Cal
Y(\Omega_2))$, where $\Omega$ is the disjoint union of $\Omega_1$
and $\Omega_2$. In particular, all torsion pairs $(\Cal X,\Cal Y)$
with $\bold t \subset \Cal X$ and $\bold q \subset \Cal Y$ split.}
    \medskip
Proof: Let $(\Cal X,\Cal Y)$ be a torsion pair with $\bold t
\subset \Cal X$ and $\bold q \subset \Cal Y$. In case the generic
module $G$ belongs to $\Cal Y$, all the Pr\"ufer modules
$S[\infty]$ belong to $\Cal Y$, since they are factor modules of
$G$, thus $\Cal Y \supseteq \Cal D,$ but then $\Cal X = \Cal R$
and $\Cal Y = \Cal D$. Now, let us assume that $G$ does not belong
to $\Cal Y$. Denote by $\Omega_2$ the set of all $\lambda\in
\Omega$ such that $\Cal T(\lambda)\cap \Cal Y$ contains a non-zero
module, and let $\Omega_1 = \Omega\setminus \Omega_2.$ If
$\lambda\in \Omega_2,$ then $\Cal Y$ contains at least one and
thus all the Pr\"ufer modules from $\Cal T(\lambda)$, since all of
them are epimorphic images of a given one. It follows that $\Cal Y
= \Cal Y(\Omega_2)$ and thus $\Cal X = \Cal X(\Omega_1).$
    \medskip
The lattice of the subcategories $\Cal Y$, where $(\Cal X,\Cal Y)$
is a torsion pair with $\bold t \subset \Cal X$ and $\bold q
\subset \Cal Y$, looks as follows, where the lower part is order
isomorphic to the power set $\Cal P(\Omega)$: 
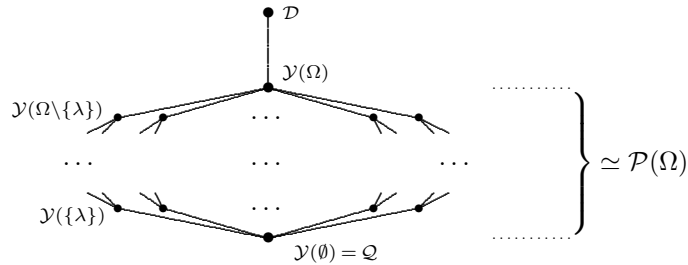
$$
\hbox{\beginpicture \setcoordinatesystem units <1cm,1cm>
%==========================================1
\put{} at 0 0 \put{} at 4 3 \put{$\bullet$} at 2 0 \put{$\ssize
\Cal Y(\emptyset)\,=\,\Cal Q $} at 2.9 -0.2 \put{$\bullet$} at 2 2
\put{$\ssize \Cal Y(\Omega)$} at 2.5 2.2 \put{$\bullet$} at 2 3
\put{$\ssize \Cal D$} at 2.3 3 \put{$\ssize \bullet$} at 0 1.6
\put{$\ssize \bullet$} at 0 0.4

\put{$\ssize \bullet$} at 0.6 1.6 \put{$\ssize \bullet$} at 0.6
0.4

\put{$\ssize \bullet$} at 3.4 1.6 \put{$\ssize \bullet$} at 3.4
0.4

\put{$\ssize \bullet$} at 4 1.6 \put{$\ssize \bullet$} at 4 0.4

\put{$\ssize \Cal Y(\{\lambda\})$} at -0.6 0.3 \put{$\ssize \Cal
Y(\Omega\setminus\{\lambda\})$} at -0.8 1.7 \plot 2 2  2 3 / \plot
2 2  0 1.6 / \plot 2 2  0.6 1.6 / \plot 2 2  3.4 1.6 / \plot 2 2
4 1.6 / \put{$\dots$} at 2 1.6 \plot 2 0  0   0.4 / \plot 2 0  0.6
0.4 / \plot 2 0  3.4 0.4 / \plot 2 0  4   0.4 /

\plot 0 1.6  -0.2 1.4 / \plot 0 1.6  -0.4 1.4 /

\plot 0 0.4  -0.2 0.6 / \plot 0 0.4  -0.4 0.6 /

\plot 4 1.6   4.2 1.4 / \plot 4 1.6   4.4 1.4 /

\plot 4 0.4   4.2 0.6 / \plot 4 0.4   4.4 0.6 /

\plot 0.6 1.6  0.5 1.4 / \plot 0.6 1.6  0.3 1.4 /

\plot 0.6 0.4  0.5 0.6 / \plot 0.6 0.4  0.3 0.6 /

\plot 3.4 1.6  3.5 1.4 / \plot 3.4 1.6  3.7 1.4 /

\plot 3.4 0.4  3.5 0.6 / \plot 3.4 0.4  3.7 0.6 /

\put{$\dots$} at 2 0.4 \put{$\dots$} at -0.5 1 \put{$\dots$} at 2
1 \put{$\dots$} at 4.5 1 \put{$\left.\phantom {\matrix a \cr b\cr
c\cr d \cr e  \endmatrix} \right\}$} at 6 1 \put{$\simeq \Cal
P(\Omega)$} at 7 1 \setdots <1mm> \plot 5 0  6 0 / \plot 5 2  6 2
/
\endpicture}
$$
    \bigskip
Let us add the following observations which will be useful in
section 10.
        \medskip
{\bf Lemma 8.} {\it Let $Y$ be in $\omega$. Any monomorphism $X
\to Y$ with $X\in \omega_0$ splits. Any epimorphism $Y \to Z$ with
$Z\in \Add G$ splits.}
        \medskip
Proof: Decompose $Y = Y_1 \oplus Y_0$ with $Y_1\in \Add G$ and
$Y_0 \in \omega_0.$ Now $\Hom(X, Y_1) = 0$, thus any map $X \to Y$
maps into $Y_0$. But clearly any monomorphism in $\omega_0$
splits. Similarly, $\Hom(Y_0,Z) = 0,$ thus any map $Y \to Z$
vanishes on $Y_0$ and induces an epimorphism $Y_1 \to Z$. But any
epimorphism in $\Add G$ splits. This completes the proof.
        \medskip
{\bf Lemma 9.} {\it The subcategories $\Cal C, \Cal D$ and
$\omega$ are closed under products.}
        \medskip
Proof: Since $\omega = \Cal C \cap \Cal D$, we only have to
consider $\Cal C$ and $\Cal D$. It is clear that $\Cal C$ is
closed under products, since $\Cal C$ is a torsionfree class. For
$\Cal D,$ we use the description $\Cal D = \{M\mid \Ext^1(\bold
t,M) = 0\}$ in order to see it.
        \bigskip\bigskip
%============================
{\bf 9\. Concealed canonical algebras.}
    \bigskip
Let us outline why the results presented above for the canonical
algebras immediately extend to the more general class of concealed
canonical algebras using tilting. As we have mentioned above,
according to [LP] this takes care of any sincere stable separating
tubular family, thus in particular of all the stable separating
tubular families of a tubular algebra.

Recall that the concealed canonical algebras are obtained in the
following way: we start with a canonical algebra $\Lambda$ with
the canonical trisection $(\bold p,\bold t, \bold q)$ and we
consider the subcategories $\Cal C, \Cal D, \omega$ as defined
above. Let $T$ be a tilting module which belongs to $\bold p$ and
consider $\Lambda' = \End(T)^{\text{op}}.$ Then, by definition,
$\Lambda'$ is a concealed canonical algebra. The tilting functor
$F = \Hom(T,-)$ sends $\bold t$ to a sincere stable separating
tubular family $\bold t'$ in $\mod \Lambda'$, and all sincere
stable separating tubular families are obtained in this way. The
tubular family $\bold t'$ separates say $\bold p'$ from $\bold
q'$, here $\bold p'$ is the image of $\bold p$ under $F$, whereas
the modules $M'$ in $\bold q'$ are extensions  of the form $$
 0 \to F'(M_1) \to M' \to F(M_2) \to 0
$$ with $M_1\in \add\bold p$, $M_2\in \add\bold q$, and $F' =
\Ext^1(T,-)$.

As above, we may define $$
 \Cal C' = r(\bold q'),\ \Cal Q' = g(\bold q'), \t{and}
 \Cal R' = rl(\bold t'),\ \Cal D' = l(\bold t').
$$ Of course, we also put $\omega' = \Cal C' \cap \Cal D'$. The
generic module $G'$ and the Pr\"ufer modules for $\Lambda'$ are
defined as the images of the corresponding $\Lambda$-modules under
$F$. But we may define the Pr\"ufer modules for $\Lambda'$ also
directly as unions of chains of indecomposable modules in $\bold
t'$. And we denote by $\omega'_0$ the full subcategory of all
direct sums of Pr\"ufer modules.

We claim that the following assertions hold:
\item{(1)} {\it The categories $\omega$ and
$\omega'$ are equivalent categories,} an equivalence being given
by the restriction of $F$, thus any object in $\omega'$ is a
direct sum of indecomposable objects and the indecomposables in
$\omega'$ are a generic module and Pr\"ufer modules. And this
equivalence yields an equivalence of $\omega_0$ and $\omega'_0$.
\item{(2)} {\it Any $\Lambda'$-module is a
direct sum of a module in $\Cal R',$ a module in $\omega'$ and a
module in $\Cal Q'$.}
\item{(3)} {\it For each module $C'$ in $\Cal C'$ there is an exact sequence
$$
  0 @>>> C' @>f>> X' @>>> Y' @>>> 0
$$ with $Y'\in \omega'_0$ and $X' \in \omega'$, where $f\:C' \to
X'$ is a minimal left $\omega$-approximation.}
\item{(4)}  {\it For each module $D'$ in $\Cal D'$ there is an exact sequence
$$
  0 @>>> X' @>>> Y' @>g>> D' @>>> 0
$$ with $X'$ a direct sum of copies of $G'$ and $Y' \in \omega'$,
where $g\:Y' \to D'$ is a minimal right $\omega$-approximation.}
\item{(5)} $\Ext^1(\Cal C',\Cal D') = 0.$
\item{(6)} {\it
$\Cal C' = \{C'\mid \Ext^1(C',\omega') = 0\}$  and $\Cal D' =
\{D'\mid \Ext^1(\omega',D') = 0\}$  and $\omega' = \{B'\in \Cal
C'\mid \Ext^1(\Cal C',B') = 0\} =
  \{A'\in \Cal D'\mid \Ext^1(A',\Cal D') = 0\}$.}
\item{(7)} {\it $\pd C' \le 1$ for any $C' \in \Cal C'$
and $\id D' \le 1$ for any $D' \in \Cal D'$.}
\item{(8)} {\it Let $Y'$ be in $\omega'$.
Any monomorphism $X' \to Y'$ with $X'\in \omega'_0$ splits. Any
epimorphism $Y' \to Z'$ with $Z'\in \Add G'$ splits.}
    \medskip
    Proof: We
recall that we denote by $r(T)$ the full subcategory of $\Mod
\Lambda$ given by the modules $M$ with $\Hom(T,M) = 0$, thus
$(r(T),g(T))$ is the torsion pair in $\Mod \Lambda$ attached to
the tilting module $T$. As usual, we also need the cotilting
module $T' = F(D\Lambda)$ in $\mod \Lambda'$ and the full
subcategory $l(T')$ of all $\Lambda'$-modules $N$ with $\Hom(N,T')
= 0,$ so that $(c(T'),l(T'))$ is the torsion pair attached to the
cotilting module $T'$. Tilting theory asserts that $F$ yields an
equivalence between $g(T)$ and $c(T')$ and that $F'$ yields an
equivalence between $r(T)$ and $l(T'),$ see [CF].

Let us look at the eight assertions. Assertion (1) is obvious,
since $\omega$ is contained in $g(T)$. In order to verify (2), the
essential observation is the following vanishing result $$
 \Ext^1(F(\Cal C), F'(r(T))) = 0.
$$ This formula may be shown directly, but one also may prefer to
work in the bounded derived category $D^b(\Mod\Lambda)$, see [R3].
It follows from the formula that any $\Lambda'$-module is the
direct sum of a module in $F(\Cal C)$ and a module $N$ which has a
submodule $ N'$ in $F'(r(T))$ such that $ N/N'$ belongs to $
F(\Cal Q)$. Then it is easy to see that $\Hom(N',\bold t') = 0 =
\Hom(N/N',\bold t')$, so that $N'$ and $N/N'$ are in $\Cal D'$,
but clearly not in $\omega'$. Since $\bold q'$ is numerically
determined, $(\Cal C',\Cal Q')$ is a split torsion pair and hence
$N$ is in $\Cal Q'$. It follows that $F(\Cal C) = \Cal C'$ by
Proposition 1. It remains to show that $F(\Cal R) \subseteq \Cal
R'$, or, in other words $\Hom(\bold q',F(\Cal R)) = 0.$ Now
$\Hom(F(\bold q),F(\Cal R)) = 0$, since $\Hom(\bold q, \Cal R) =
0$ and also $\Hom(F'(M_1),F(\Cal R)) = 0$, for $M_1\in  \add \bold
p.$

For (3) we use that $F\:\Mod\Lambda \to \Mod\Lambda'$ restricts to
an equivalence from $\Cal C\cap g(T)$ to $\Cal C'$, and from
$\omega$ to $\omega'$, together with the corresponding result for
$\Mod\Lambda$.

The proof of Theorem 5 generalizes to $\Lambda'$-modules. We need
only to observe that $\Ext^1(\omega',\omega') = 0.$ This follows
from $\Ext^1(\omega,\omega) = 0$ by using the inverse equivalence
from $\omega'$ to $\omega$ induced by $T\otimes_{\Lambda'}\: \Mod
\Lambda' \to \Mod\Lambda$. Hence we get (4).

We have already pointed out that the torsion pair $(\Cal C',\Cal
Q')$ splits, and it follows from (2) that the pair $(\Cal R',\Cal
D')$ splits. Hence we have $\Ext^1(\Cal C',\Cal Q') = 0 =
\Ext^1(\Cal R',\Cal D').$ We have also seen that
$\Ext^1(\omega',\omega') = 0.$

We claim that the modules in $\omega'$ have projective and
injective dimension at most one. The Pr\"ufer modules have
projective and injective dimension one, since they are direct
limits of modules in $\bold t'$ which have projective and
injective dimension one. Applying the exact sequences in (3) and
(4) to $R' = \Lambda$ and $Q' = D\Lambda$, we see that also the
generic module $G'$ has projective and injective dimension one.
Hence the functors $\Ext^1(M',-)$ and $\Ext^1(-,M')$ are right
exact for $M' \in \omega'$. Since $\Cal C' = c(\omega')$ by (3)
and $\Cal D' = g(\omega')$ by (4), it follows that $\Ext^1(\Cal
C',\omega') = 0 = \Ext^1(\omega',\Cal D').$ This proves (5), and
(6) now follows easily (see Propositions 2 and 3). For (7), see
the proof of Corollary 2, for (8) that of Lemma 8.
    \bigskip\bigskip
%======================================================================
{\bf 10\. Inf-tilting  and inf-cotilting modules.}
        \medskip
In this section we discuss connections with tilting theory, for
concealed canonical algebras. A usual (finite length) tilting
module yields a torsion pair, but not all torsion pairs are
obtained in this way. We are going to show that in our situation
some generalization of the concept of a tilting module which
allows a tilting  module to be of infinite length is very helpful.
In order to distinguish this generalization from the traditional
notion we refer to these modules as ``inf-tilting''  and
``inf-cotilting''  modules.
        \medskip
Up to now, the torsion pairs which we have considered explicitly
were torsion pairs in a complete module category $\Mod R$, where
$R$ is any ring. Of course, implicitly, we also dealt with torsion
pairs in categories of the form $\mod \Lambda$ with $\Lambda$ an
artin algebra. Indeed, the general concept of a torsion pair
$(\Cal F,\Cal G)$ is defined in an arbitrary abelian category
$\Cal A$; one requires that $\Hom(G,F) = 0$ for all $F\in \Cal F$
and $G\in \Cal G$, that $\Cal F$ and $\Cal G$ are closed under
isomorphisms and that for every object $A\in \Cal A$ there exists
a short exact sequence $0 \to A' \to A \to A'' \to 0$ with $A' \in
\Cal G$ and $A''\in \Cal F$. From now on, we will use the
operators $r(-),\; l(-),\; g(-),\; c(-)$ in this more general
setting and we hope that this will not lead to any confusion.
        \medskip
Given an artin algebra $\Lambda$, torsion pairs in the category
$\mod\Lambda$ of modules of finite length occur frequently as
being related to a tilting or a cotilting module. Given a tilting
module $T$ of projective dimension at most one, the pair
$(r(T),g(T))$ in $\mod\Lambda$ is a torsion pair. We say that it
is {\it associated} with the tilting module $T$. Similarly, given
a cotilting module $T$ of injective dimension at most one, the
pair $(c(T),l(T))$ in $\mod\Lambda$ is a torsion pair. We say that
it is {\it associated} with the cotilting module $T$. Starting
with tilting or cotilting modules, we obtain in this way many
torsion pairs, but there are also interesting torsion pairs $(\Cal
Y, \Cal X)$ in $\mod\Lambda$ which do not appear in this way.
Especially interesting are those where $\Cal X$ contains the
injective modules or $\Cal Y$ contains the projective modules. In
[HRS], these are called tilting and cotilting torsion pairs
respectively, and it is possible to imitate the usual tilting
procedure passing from $\Lambda$ to $\Gamma=\End(T)^{\text{op}} $
by performing tilting  with respect to such a torsion pair inside
the bounded derived category.

An example of a tilting and cotilting torsion pair in
$\mod\Lambda$ not associated with a tilting or cotilting
$\Lambda$-module is $(\add(\bold p\cup\bold t),\,\,\add\bold q)$,
where $\Lambda$ is a canonical algebra with canonical trisection
$(\bold p,\bold t,\bold q)$ (or more generally any concealed
canonical algebra \dots). Let us turn our attention to arbitrary,
not necessarily finitely generated modules, and consider the two
extremal torsion pairs $(\Cal R,\Cal D)$ and $(\Cal C,\Cal Q )$
which are extensions to $\Mod\Lambda$ of the torsion pair
$(\add(\bold p\cup\bold t),\,\,\add\bold q)$ in $\mod\Lambda$. We
claim that these torsion pairs are associated to something like
tilting and cotilting modules respectively: we need to work with a
generalization of the concept of a tilting or a cotilting module
which allows to deal with infinitely generated modules. Let us
refer here to Colpi-Trlifaj [CP] where these inf-tilting  modules
of projective dimension at most 1 have been considered and to
Colpi-D'Este-Tonolo [CET] for an investigation of inf-cotilting
modules of injective dimension at most 1, but also to [AC] and
[ATT].

{\bf Definition:} Let $R$ be any ring. An $R$-module $W$ of
projective dimension at most one will be called an {\it
inf-tilting} module, provided it satisfies the following two
properties: We have $\Ext^1(W,\Add W)=0$ and there is an exact
sequence $0\to \Lambda\to X\to Y\to 0$ with $X$ and $Y$ in $\Add
W$. Dually, an $R$-module $W$ of injective dimension at most one
will be called an {\it inf-cotilting} module provided
$\Ext^1(\Prod W,W) = 0$ and there is an exact sequence $0\to X'\to
Y'\to D\Lambda\to 0$ with $X'$ and $Y'$ in $\Prod W$. Here,
$D\Lambda$ denotes the dual of $\Lambda_\Lambda$, this is the
minimal injective cogenerator (at least in case $\Lambda$ is
basic).

Now assume again that $\Lambda$ is an artin algebra and let $\bold
t$ be a sincere stable tubular family in $\mod \Lambda$ separating
$\bold p$ from $\bold q$ (in particular, $\Lambda$ is concealed
canonical). Let $\Cal C = r(\bold q)$ and $\Cal D = l(\bold t)$.
The crucial subcategory to be considered is $\omega=\Cal C\cap\Cal
D$. Let $W_0$ be the direct sum of all the Pr\"ufer modules in
$\omega$, one copy from each isomorphism class, and let $W =
G\oplus W_0$, where $G$ is the generic module in $\omega$.
        \medskip
{\bf Proposition 5}. {\it The module $W$ is an inf-tilting module
of projective dimension one, the modules $W$ and $W_0$ are
inf-cotilting modules of injective dimension one and } $$
 \omega = \Add W = \Prod W = \Prod W_0.
$$
        \medskip
Proof: The following references are all to section 9. The fact
that $W$ and $W_0$ have projective dimension one and injective
dimension one follows from (7), since both modules belong to
$\omega=\Cal C\cap\Cal D$ (and are neither projective nor
injective).

Clearly, $\omega = \Add W$. According to Lemma 9, $\omega$ is
closed under products, thus we have $\Prod W_0 \subseteq \Prod W
\subseteq \omega$. In order to show $\omega \subseteq \Prod W_0$,
it is sufficient to verify that any direct sum $\bigoplus_I W$
with $I$ an infinite index set belongs to $\Prod W_0.$ Thus, let
$I$ be an infinite index set and consider the product $Y = \prod_I
W_0$ of $I$ copies of $W_0$. We know that $Y$ belongs to $\omega$,
thus it is a direct sum of copies of the generic module and the
Pr\"ufer modules. It is sufficient to show that in such a direct
sum decomposition, any indecomposable module occurs with
multiplicity at least $I$. First, consider a Pr\"ufer module $P$.
We have obvious inclusion maps $\bigoplus_I P \to \prod_I P \to
\prod_I W_0 = Y.$ According to (8), this monomorphism splits.
Also, as Krause ([K], see also [R6]) has shown, $\prod_I P$
contains $\bigoplus_I G$ as a direct summand, thus we obtain an
epimorphism $Y = \prod_I W_0 \to \prod_I P \to \bigoplus_I G.$ We
use the second part of (8) in order to conclude that $\bigoplus_I G$
is a direct summand of $Y$. Altogether we see that $\bigoplus_I W$
is a direct summand of $Y$.

The remaining assertions now follow easily: Since $\Add W = \omega
= \Prod W$, it follows from the basic splitting theorem that
$\Ext^1(W,\Add W) = 0 = \Ext^1(\Prod W,W)$ and of course also
$\Ext^1(\Prod W_0,W_0) = 0.$ Since ${}_\Lambda\Lambda$ belongs to
$\bold p \subseteq \Cal C$, its minimal left
$\omega$-approximation yields an exact sequence $0\to \Lambda\to
X\to Y\to 0$ with $X$ and $Y$ in $\omega = \Add W$, see (3). This
shows that $W$ is an inf-tilting module of projective dimension
one. Dually, the module $D\Lambda$ belongs to $\bold q \subseteq
\Cal D$, thus its minimal right $\omega$-approximation yields an
exact sequence $0\to X'\to Y'\to D\Lambda\to 0$ with $X'$ and $Y'$
in $\omega =\Prod W_0 = \Prod W$, see (4). Thus both $W$ and $W_0$
are inf-cotilting modules of injective dimension one.
        \medskip
Note that the torsion pair $(\Cal R,\,\Cal D)$ is associated with
the inf-tilting  module $W$, since $\Cal D=g(W)$, and ($\Cal
C,\Cal Q )$ is associated with the inf-cotilting  modules $W$ and
$W_0$, since $\Cal C=c(W) = c(W_0)$.
        \bigskip
%==========================================
{\bf Remark.} Note that the torsion pair $(\Cal R,\Cal D)$ does
not seem to be associated with something like a cotilting module,
but all the torsion pairs $(\Cal X(\complement\Omega'), \Cal
Y(\Omega'))$ are, where $\Omega'$ is a subset of $\Omega$ and
$\complement\Omega'$ is its complement inside $\Omega$. Namely,
{\it define $T(\Omega')$ as the direct sum of the generic module
$G$, the Pr\"ufer modules $S[\infty]$ with $S\in \Cal T(\Omega')$
and the adic modules $\widehat S$ with $S \in \Cal T(\complement
\Omega')$. Then $T(\Omega')$ is an inf-cotilting  module and $\Cal
X(\complement\Omega') = c(T(\Omega')).$} (The adic module
$\widehat S$ is the inverse limit of a chain of epimorphisms $$
 \cdots \to [n]S \to [n-1]S \to \cdots \to [2]S \to S,
$$ where $[n]S$ is the (uniquely determined) module in $\bold t$
of regular length $n$ which has $S$ as a factor module, see for
example [R5].) In the case of a tame hereditary algebra, we may
refer to [BK] for a description of all the pure injective
cotilting modules.
    \bigskip\bigskip
%=========================
{\bf 11\. Derived equivalent categories.}
        \medskip
In this section we outline the effect of tilting with respect to
some of the torsion pairs considered above inside the derived
category $D^b(\Mod\Lambda).$
        \medskip
If $R$ is any ring, let  $D^b(\Mod R)$ be its bounded derived
category (with shift automorphisms $X \mapsto X[n]$ for all $n\in
\Bbb Z$ and homology functors $H^n\:D^b(\Mod R) \to \Mod R$). We
always will identify $\Mod R$ with the full subcategory of all
objects $X$ in $D^b(\Mod R)$ with $H^i(X) = 0$ for $i\neq 0.$
Given a torsion pair $(\Cal X,\Cal Y)$ in $\Mod R$, there is an
inf-tilting  procedure inside $D^b(\Mod R)$ with respect to this
torsion pair. It yields a new abelian category $\Cal A$ which is
contained in $D^b(\Mod R)$, as follows: $\Cal A = \Cal A(\Cal
X,\Cal Y)$ is the full subcategory of all objects $A$ of $D^b(\Mod
R)$ such that $$
 H^{-1}(A) \in \Cal Y,\quad H^0(A)\in \Cal X \t{and} H^i(A) = 0 \ \text{for}\ i\notin
\{-1,0\} $$ (see [HRS]). Under the condition that $\Cal X$
contains all projective $R$-modules or that $\Cal Y$ contains all
injective $R$-modules, it follows that $\Cal A$ is derived
equivalent to $\Mod R$ and that $(\Cal Y[-1],\Cal X)$ is a torsion
pair in $\Cal A$. In case the torsion pair $(\Cal X,\Cal Y)$ is
split with $\pd X \le 1$ for $X$ in $\Cal X$, and if $\Cal X$
contains all the projectives, then the new abelian category is
hereditary (see [HR] and [HRS]).
    \medskip
Let us now assume again that $\Lambda$ is a canonical algebra with
a stable tubular family $\bold t$ separating $\bold p$ from $\bold
q$. We consider the subcategories $\Cal C, \Cal D, \omega,$ and so
on, as defined above, relative to $\bold t$. In particular, let us
look at some of the torsion pairs $(\Cal X,\Cal Y)$ with $\bold
t\subset \Cal X$ and $\bold q \subset \Cal Y$. It is interesting
to observe that some of the exact sequences in $\Mod \Lambda$
which have been discussed in this paper, can be interpreted as
injective or projective resolutions in the new abelian category
$\Cal A$, depending on the choice of the torsion pair $(\Cal
X,\Cal Y)$ and that the generic module $G$ and the Pr\"ufer
modules yield enough injective or projective objects in $\Cal A$.
    \medskip
{\bf Proposition 6.} {\it The category $\Cal A = \Cal A(\Cal C,
\Cal Q )$ is a hereditary abelian category derived equivalent to
$\Mod \Lambda$. The pair $(\Cal Q [-1],\Cal C)$ is a torsion pair
in $\Cal A$. The subcategory $\omega$ is the class of all
injective objects in $\Cal A$ and $\Cal A$ has sufficiently many
injective objects.}
    \medskip
{\bf Proposition 7.} {\it The category $\Cal A' = \Cal A(\Cal
R,\Cal D)$ is a hereditary abelian category derived equivalent to
$\Mod \Lambda$. The pair $(\Cal D[-1],\Cal R)$ is a torsion pair
in $\Cal A'$. The subcategory $\omega[-1]$ is the class of all
projective objects in $\Cal A'$ and $\Cal A'$ has sufficiently
many projective objects.}
    \medskip
{\bf Proposition 8.} {\it The category $\Cal A'' = \Cal A(\Cal
X(\emptyset),\Cal Y(\Omega))$ is a hereditary abelian category
derived equivalent to $\Mod \Lambda$. The pair $(\Cal
Y(\Omega)[-1],\Cal X(\emptyset))$ is a torsion pair in $\Cal A''$.
If $P$ is a Pr\"ufer module in $\Mod\Lambda$, then $P[-1]$ is an
indecomposable projective object of $\Cal A$, whereas the generic
module $G$, considered as an object in $\Cal A$ is simple
injective. Thus, $\Cal A''$ has non-zero projective and non-zero
injective objects, but neither sufficiently many projective
objects nor sufficiently many injective objects.}
    \medskip
The proofs of these propositions follow quite easily from the
above remarks and the properties of torsion pairs in question
which have been established in previous sections.

Note that in the opposite direction it is shown in [L] that
generic modules over canonical algebras can be investigated by
first considering generic sheaves in the category of quasicoherent
sheaves over weighted projective lines.
    \bigskip\bigskip

{\bf 12\.} {\bf Additional comments}.
        \medskip
In this section we discuss the relationship between split torsion
pairs of $\mod\Lambda$ and $\Mod\Lambda$. We also indicate briefly
a different approach to the study of $\Mod\Lambda$ when $\Lambda$
is a canonical algebra, using tame bimodules.

In this paper we have considered in detail the cut of $\mod
\Lambda$ between $\bold t$ and $\bold q$, where $\bold t$ is a
sincere stable tubular family separating $\bold p$ from $\bold q$:
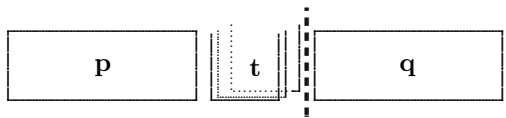
$$ \hbox{\beginpicture \setcoordinatesystem units <.9cm,.4cm>
\put{} at 0 0.5 \put{} at 6 3 \plot 0 0   2.78 0  2.78 2.3    0
2.3  0 0 / \plot 4.53 0 7.3 0  7.3 2.3  4.53 2.3  4.53 0 / \plot 3
2.2  3 0  4 0  4 2.2 / \plot 4.1 0.1  4.1 2.3 / \plot 4.2 0.3  4.3
0.3  4.3 2.5 / \setdots <1pt> \plot 3.1 2.3 3.1 0.1  4.1 0.1 /
\setdots <2pt> \plot 3.3 2.5 3.3 0.3  4.2 0.3 / \put{$\bold p$} at
1.4 1.1 \put{$\bold t$} at 3.65 1.1 \put{$\bold q$} at 5.9 1.1
\setdashes<1mm> \plot 4.41 -0.5 4.41 3.1 / \plot 4.43 -0.5 4.43
3.1 / \plot 4.39 -0.5 4.39 3.1 /
\endpicture}
$$
 One may try to look also at the dual cut between $\bold p$
and $\bold t$: 
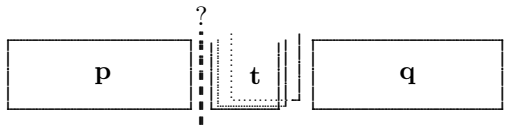
$$ \hbox{\beginpicture \setcoordinatesystem units
<.9cm,.4cm> \put{} at 0 -.5 \put{} at 6 3 \plot 0 0   2.7 0  2.7
2.3    0 2.3  0 0 / \plot 4.5 0 7.3 0  7.3 2.3  4.5 2.3  4.5 0 /
\plot 3 2.2  3 0  4 0  4 2.2 / \plot 4.1 0.1  4.1 2.3 / \plot 4.2
0.3  4.3 0.3  4.3 2.5 / \setdots <1pt> \plot 3.1 2.3 3.1 0.1  4.1
0.1 / \setdots <2pt> \plot 3.3 2.5 3.3 0.3  4.2 0.3 / \put{$\bold
p$} at 1.4 1.1 \put{$\bold t$} at 3.65 1.1 \put{$\bold q$} at 5.9
1.1 \setdashes<1mm> \plot 2.85 -0.5 2.85 2.6 / \plot 2.87 -0.5
2.87 2.6 / \plot 2.83 -0.5 2.83 2.6 / \put{?} at 2.85 3.1
\endpicture}
$$

 \bigskip
The situation seems to be similar, but it is not! There is a dual
cut only when dealing with finite dimensional representations ---
the behavior of the infinite dimensional modules in this part of
the category $\Mod\Lambda$ is far more complicated: indeed, there
do exist many torsion pairs $(\Cal X,\Cal Y)$ in $\Mod\Lambda$
with $\bold p \subset \Cal X$ and $\bold t \subset \Cal Y$ which
do not split (indeed, we do not know any one which splits).

In order to provide at least one example, let us consider again
the special case where $\Lambda$ is the Kronecker algebra, thus we
consider the representations of the quiver $$ \hbox{\beginpicture
\setcoordinatesystem units <0.8cm,0.5cm>
%==========================================1
\put{$\circ$} at 0 0 \put{$\circ$} at 2 0 \put{$a$} at -0.3 0
\put{$b$} at 2.3 0 \put{$\ssize\alpha$} at 1 0.7
\put{$\ssize\beta$} at 1 -0.8
%======================
\arrow <2mm> [0.25,0.75] from 0.6 -0.3 to 0.2 -0.2 \arrow <2mm>
[0.25,0.75] from 0.6  0.35 to 0.2  0.25
%\arrow <2mm> [0.25,0.75] from 1.8  0 to 0.2  0
%======================
\setquadratic \plot 1.8 -0.2   1.2 -0.35  0.6  -0.3  / \plot 1.8
0.25   1.2  0.4  0.6 0.35  /
%======================
\endpicture}
$$ We consider the torsion pair $(\Cal X,\Cal Y)$, where $\Cal X =
r(\bold t)$  and $\Cal Y = g(\bold t).$ Note that the full
subcategory $\Cal M$ of all representations $V$ with $V(\beta)$
being an identity map is isomorphic to the category of
$k[T]$-modules (an isomorphism is obtained by sending the
$k[T]$-module $M$ to the representation $V_M$ with $V_M(a) =
V_M(b) = M,$ such that $V_M(\alpha)$ is the multiplication by $T$
and $V_M(\beta) = 1_M$). It is quite obvious that the restriction
of the torsion pair $(\Cal X,\Cal Y)$ to $\Cal M$ is just the
usual pair of torsionfree and torsion $k[T]$-modules and it is
well-known that there do exist many $k[T]$-modules whose torsion
submodule does not split off (see [F], Chapter XIV).

Given such a trisection $(\bold p,\bold t,\bold q)$ of $\mod \Lambda$, the
difference of the two cuts between $\bold p$ and $\bold t$ on one
hand and between $\bold t$ and $\bold q$ on the other hand should
not prevent a detailed study of what lies inbetween $\bold p$ and
$\bold q$ in $\Mod\Lambda$, namely the subcategory $\Cal M(\bold t)=l(\bold p)\cap r(\bold
q)$. For the case of a tubular algebra $\Lambda$ and its various
trisections $(\bold p_w,\bold t_w,\bold q_w)$, with $w\in
Q^\infty_0$, the subcategory $\Cal M(\bold t_w)$ can be denoted
just by $\Cal M(w)$. This subcategory, as well as corresponding
 subcategories $\Cal M(w)$ for $w\in \Bbb{R}^\infty_0\backslash
Q^\infty_0$ will be studied in the final section of the paper.

Since most of the results in this paper have been shown in [R1]
for tame hereditary algebras, it would seem reasonable to
establish the general case by reducing the investigation of the
module theory for canonical algebras to that of tame bimodules,
where the corresponding assertions are already known. This
definitely can be done. Let $P_0 = S$ be the simple projective
module and $P_\infty$  the projective cover of the simple
injective module $S'$. If $P = P_0\oplus P_\infty$, then basic
properties of the functor
$F=\Hom(P,\,\,):\Mod\Lambda\to\Mod\Lambda_0$ for
$\Lambda_1=\End(P)^{\text{op}}$ were investigated in [R4] for
finitely generated modules and can be generalized without problems
to arbitrary modules. One considers the modules $C$ in
$\Mod\Lambda$ with $F(C')\neq0$ for each nonzero summand $C'$ of
$C$, and one can compare crucial properties for $C$ and $F(C)$ in
this case. This allows to use results for $\Mod\Lambda_1$ and
transfer them to $\Mod\Lambda$. Note that $\Lambda_1$ is given by
a tame bimodule, namely by $\Hom(P_0,P_\infty),$ so that all the
relevant properties of the category $\Mod\Lambda_1$ are known for
a long time. In particular, in this way the canonical generic
module does not have to be constructed from scratch for all the
canonical algebras, but only for the tame bimodules involved. On
the other hand, we hope that the direct approach presented in this
paper helps to trace in which way the structure of the category of
finite dimensional representations determines that of all the
modules.
    \bigskip\bigskip
%======================================================================
{\bf 13\. Tubular algebras}
        \medskip
   \medskip
Let $\Lambda$ be a canonical algebra with canonical trisection
$(\bold p, \bold t,\bold q)$. Let $S$ be the simple projective
module and $S'$ the simple injective module. We denote by
$\Lambda_0$ the factor algebra of $\Lambda$ so that
$\mod\Lambda_0$ is the full subcategory of all $\Lambda$-modules
$M$ with $[M:S']= 0.$ Similarly, $\Lambda_\infty$ is the factor
algebra of
 $\Lambda$ with the property that
$\mod\Lambda_\infty$ is the full subcategory of all
$\Lambda$-modules $M$ with $[M:S]= 0.$ Note that both $\Lambda_0$
and $\Lambda_\infty$ are hereditary algebras. The representation
types of $\Lambda_0$ and of $\Lambda_\infty$ coincide and
determine the representation type of $\Lambda$. Let us review the
different cases.

If $\Lambda_0$ and $\Lambda_\infty$ are of finite representation
type, then $\Lambda$ is a tame concealed algebra, $\bold p$ is a
preprojective component, $\bold q $ a preinjective component. This
is essentially the case which has been studied in detail in [R1].
In particular, it has been shown there that $\Cal Q = \Add \bold
q$, whereas $\Cal R$ is a wild category (see also [R7]). In this
case, the asymmetry between $\Cal R$ and $\Cal Q$ is best visible.

If $\Lambda_0$ and $\Lambda_\infty$ are of wild representation
type, then $\Lambda$ is of course also of wild representation
type. In this case, not much is known even for the finite
dimensional $\Lambda$-modules, but it should be worthwhile to
study this case in more detail in future.

It remains to consider the case where both $\Lambda_0$ and
$\Lambda_\infty$ are of tame representation type, in this case
$\Lambda$ is said to a {\it tubular canonical} algebra. More
generally, we may consider an arbitrary {\it tubular} algebra,
these are the concealed canonical algebras obtained from a tubular
canonical algebra by tilting.
        \medskip
>From now on, let $\Lambda$ be a tubular algebra. The structure of
$\mod \Lambda$ is known in detail (see [R2] and [LP]). There is a
preprojective component $\bold p_0$ and a preinjective component
$\bold q_\infty$. We denote by $I_0$ the ideal which is maximal
with the property that it annihilates all the modules in $\bold
p_0$ and by $I_\infty$ the ideal which is maximal with the
property that it annihilates all the modules in $\bold q_\infty$.
Then we obtain factor algebras $\Lambda_0 = \Lambda/I_0$ and
$\Lambda_\infty = \Lambda/I_\infty$ which both are tame concealed
algebras (in the case of a tubular canonical algebra, we recover
the factor algebras already introduced). Let $\bold t_0$ be the
Auslander-Reiten components of $\mod\Lambda$ which contain regular
$\Lambda_0$-modules, and $\bold t_\infty$ those which contain
regular $\Lambda_\infty$-modules. Then both $\bold t_0$ and $\bold
t_\infty$ are sincere separating tubular families, but both are
not stable ($\bold t_0$ will contain indecomposable projective
modules, $\bold t_\infty$ indecomposable injective ones). If we
denote by $\bold q_0$ the indecomposable modules in $\mod\Lambda$
which do not belong to $\bold p_0$ or $\bold t_0$, then $\bold
t_0$ separates $\bold p_0$ from $\bold q_0$. If we denote by
$\bold p_\infty$ the indecomposable modules in $\mod\Lambda$ which
do not belong to $\bold t_\infty$ or $\bold q_\infty$, then $\bold
t_\infty$ separates $\bold p_\infty$ from $\bold q_\infty$. The
modules in $\bold q_0\cap \bold p_\infty$ fall into a countable
number of sincere stable separating tubular families $\bold
t_\alpha$, indexed by $\alpha\in \Bbb Q^+$, such that for $\alpha
< \beta$ in $\Bbb Q^+$ the class $\bold t_\alpha$ generates $\bold
t_\beta$, and also $\bold t_\alpha$ is cogenerated by $\bold
t_\beta$. More generally, this generation and cogeneration
property holds for all $\alpha < \beta$ in $\Bbb Q_0^\infty =
 \Bbb Q^+\cup \{0,\infty\}$.

Let $\Bbb R_0^\infty = \Bbb R^+\cup \{0,\infty\}$. For any $w \in
\Bbb R_0^\infty$, we denote by $\bold p_w$ the modules which
belong to $\bold p_0$ or to some $\bold t_\alpha$ with $\alpha <
w$, and we denote by $\bold q_w$ the modules which belong to
$\bold t_\gamma$ with $w < \gamma$ or to $\bold q_\infty$ (here,
$\alpha, \gamma$ belong to $\Bbb Q_0^\infty$). For $\beta\in \Bbb
Q_0^\infty$ we obtain in this way a trisection $(\bold
p_\beta,\bold t_\beta,\bold q_\beta)$ of $\mod \Lambda$, with
$\bold t_\beta$ a tubular family which separates $\bold p_\beta$
from $\bold q_\beta$, and $\bold t_\beta$ is stable provided $0 <
\beta < \infty.$ For $w\in \Bbb R_0^\infty\setminus \Bbb
Q_0^\infty$, the two module classes $\bold p_w$ and $\bold q_w$
comprise all the indecomposables from $\mod\Lambda$.
        \medskip
Let us turn our attention now to arbitrary, not necessarily finite
dimensional modules. For any $w \in \Bbb R_0^\infty$, let $\Cal
C_w = r(\bold q_w)$ and $\Cal B_w = l(\bold p_w).$ The
subcategories we are interested in are those of the form $$
 \Cal M(w) = \Cal C_w \cap \Cal B_w = r(\bold q_w) \cap l(\bold p_w),
$$ defined for any $w\in \Bbb R_0^\infty$. The modules in $\Cal
M(w)$ are said to have {\it slope} $w$. Of course, for $\alpha\in
\Bbb Q_0^\infty$ the modules in $\bold t_\alpha$ as well as those
in $\omega_\alpha$ have slope $\alpha$. For non-rational $w$,
examples of modules in $\Cal M(w)$ will be presented at the end of
the section.
        \medskip
{\bf Theorem 6.} {\it Any indecomposable $\Lambda$-module which
does not belong to $\bold p_0$ or $\bold q_\infty$ has a slope.
For $0 \le w < w' \le \infty$, we have $\Hom(\Cal M(w'),\Cal M(w))
= 0.$}
        \medskip
Note that the second assertion immediately implies that $\Cal
M(w)\cap \Cal M(w') = 0,$ thus {\it if a module has a slope, its
slope is a well-defined element of $\Bbb R_0^\infty.$}
        \medskip
Before we start with the proof, let us analyze the two module
classes $\Cal C_w$ and $\Cal B_w$, as well as related ones.
        \medskip
{\bf The torsion pair $(\Cal C_w,\Cal Q_w)$.} First, we consider
$\Cal C_w = r(\bold q_w)$. Let $\Cal Q_w = g(\bold q_w).$ Note
that $\bold q_w$ is always closed under successors, thus $(\Cal
C_w,\Cal Q_w)$ is a torsion pair, according to Lemma 3 and Lemma
4. {\it For $\beta \in \Bbb Q_0^\infty$, the torsion pair $(\Cal
C_\beta,\Cal Q_\beta)$ is split.} Proof: For $\beta \in \Bbb Q^+$
and for $\beta = \infty$, the module class $\bold q_\beta$ is
numerically determined, thus we can use Proposition 1. The class
$\bold q_0$ is never numerically determined, but it is at least
numerically almost determined (the corresponding function $\delta$
vanishes precisely on those modules in $\bold q_0$ which do not
have any simple $\Lambda_0$-module as composition factor, but
there are only finitely many isomorphism classes of indecomposable
modules of this kind). Thus we can use the Remark at the end of
section 1.
        \medskip
{\bf The subcategories  $\omega_\beta$ for $\beta\in \Bbb Q^+.$}
For $\beta \in \Bbb Q^+$, the trisection $(\bold p_\beta,\bold
t_\beta, \bold q_\beta)$ allows us to use all the previous results
of the paper. In particular, there is a corresponding subcategory
$\omega_\beta$ containing a generic module $G_\beta$ as well as
Pr\"ufer modules. Actually, there are generic modules $G_\beta$
also for $\beta\in\{0,\infty\}$, thus for all $\beta \in \Bbb
Q_0^\infty$; namely, $G_0$ is the generic module of $\Lambda_0$,
and similarly, $G_\infty$ is the generic module of
$\Lambda_\infty$. According to Corollary 8, we have $l(\bold
t_\beta) = g(G_\beta)$, for $\beta\in \Bbb Q^+$, and this also
holds for $\beta = \infty$ (but not for $\beta = 0;$ in order to
show that $l(\bold t_\infty) = g(G_\infty)$, one first should
notice that both classes are contained in $\Mod \Lambda_\infty$
and then use Corollary 8 for the unique separating tubular family
of $\mod \Lambda_\infty$).
        \medskip
{\bf Lemma 10.} {\it Let $\alpha < \beta$ in $\Bbb Q_0^\infty.$
Then $\bold t_\alpha$ generates $G_\beta$. If in addition $0 <
\alpha$, then $G_\alpha$ generates $\bold t_\beta.$}
        \medskip
Proof. In order to show that $\bold t_\alpha$ generates $G_\beta$
for all $0 \le \alpha < \beta \le \infty$, consider first the case
$0 < \alpha$.

First, we claim that $G_\beta$ cannot belong to $\Cal C_\alpha$.
Choose $\gamma$ with $\alpha < \gamma < \beta$ and take an
indecomposable module $M$ in $\bold t_\gamma$. If $\beta <
\infty$, consider the left $\omega_\beta$-approximation $M \to
M_{\omega_\beta}$. Since $M_{\omega_\beta}$ is a non-zero direct
sum of copies of $G_\beta$, there are non-zero maps $M \to
G_\beta$. Consider now the case of $\beta = \infty$. Since
$\Lambda$ is a coray coextension of the tame concealed algebra
$\Lambda_\infty$, the trisection $(\bold p_\infty,\bold
t_\infty,\bold q_\infty)$ is obtained as follows: $\bold q_\infty$
consists of the preinjective $\Lambda_\infty$-modules, whereas
$\bold p_\infty$ consists of all those indecomposable
$\Lambda$-modules $N$ whose restriction $N^{(\infty)}$ to
$\Lambda_\infty$ is a direct sum of preprojective
$\Lambda_\infty$-modules. Of course, $N^{(\infty)}$ is the maximal
factor module of $N$ which is a $\Lambda_\infty$-module. Note that
$N^{(\infty)} = 0$ only for finitely many $\Lambda$-modules and all of
them belong to $\bold p_0$. Since the module $M$ belongs to $\bold
t_\gamma$, and $0 < \gamma < \infty$, we see that $M^{(\infty)}$ is a
non-zero preprojective $\Lambda_\infty$-module and therefore
$\Hom(M^{(\infty)},G_\infty) \neq 0$. Since $M^{(\infty)}$ is a factor module
of $M$, we conclude that $\Hom(M,G_\infty) \neq 0.$ Always,
$M\in\bold q_\alpha$, thus we see that $G_\beta$ can not belong to
$\Cal C_\alpha$.

Since $G_\beta$ cannot belong to $\Cal C_\alpha$, and $G_\beta$ is
indecomposable, it belongs to $\Cal Q_\alpha$. This shows that
there is a direct sum $\bigoplus_{i\in I} M_i$ of modules $M_i\in
\bold q_\alpha$ which maps onto $G_\beta.$ However, the projective
cover $P(M_i) \to M_i$ factors through $\add \bold t_\alpha$, thus
we see that any $M_i$ is generated by $\bold t_\alpha$. This shows
that $G_\beta$ is generated by $\bold t_\alpha$.

If $\alpha = 0,$ then choose $0 < \alpha' < \beta.$ By the
previous considerations, $G_\beta$ is generated by $\bold
t_{\alpha'}.$ Since $\bold t_{\alpha'}$ is generated by $\bold
t_0$, we conclude that $G_\beta$ is generated by $\bold t_0$.

In order to show the second assertion, note that we deal with $0 <
\alpha < \beta \le \infty.$ Now $\bold t_\beta \subset l(\bold
t_\alpha) = g(G_\alpha)$.
        \medskip
{\bf Remarks.} The first assertion of Lemma 10 can be strengthened
as follows: {\it If $\alpha < \beta$ in $\Bbb Q_0^\infty$ and
$\lambda \in \Omega_\alpha$, then the class $\bold
t_\alpha(\lambda)$ generates $G_\beta$} (here $\Omega_\alpha$ is
the index set for the tubular family $\bold t_\alpha$). This
follows from the proof, but can be derived also from the statement
itself: Let $\alpha < \alpha' < \beta.$ Then Lemma 10 asserts that
$\bold t_{\alpha'}$ generates $G_\beta$, but it is well-known that
any $\bold t_\alpha(\lambda)$ generates $\bold t_{\alpha'}.$

Also, let us stress that $G_0$ does not generate $\bold t_\beta$
for any $\beta\in \Bbb Q_0^\infty,$ since $G_0$ is a
$\Lambda_0$-module, whereas all the $\bold t_\beta$ contain
modules which are not $\Lambda_0$-modules. If we denote by $P$ the
direct sum of all indecomposable projective modules in $\bold
t_0$, then $G_0\oplus P$ generates $\bold t_\beta$, for any $0 <
\beta$.
        \bigskip
{\bf The module class $\Cal B_w$}. For any $w\in \Bbb
R_0^\infty$, we have defined $\Cal B_w = l(\bold p_w)$. By
definition, this is the torsion class of a torsion pair, the
corresponding torsionfree class is $r(\Cal B_w) = rl(\bold
p_w).$ For $w = 0$, the module class $\Cal B_w$ consists of all
the modules $M$ which do not have an indecomposable direct summand
in $\bold p_0$.
        \medskip
{\bf Lemma 11.} {\it Let $w\in \Bbb R^+\cup \{\infty\},$ then} $$
\align
 \Cal B_w & = \; \bigcap\nolimits_{v<w} \Cal Q_{v} \cr
 &=
  \{M\mid
    M \text{\ is generated by\ } t_\alpha
      \text{\ for any\ } \alpha\in  \Bbb Q
      \text{\ with\ } 0 <\alpha < w \} \cr
 &=
  \{M\mid
    M \text{\ is generated by\ } G_\alpha
      \text{\ for any\ } \alpha\in  \Bbb Q
      \text{\ with\ } 0 <\alpha < w \}, \cr
\endalign
$$ {\it here the $v$ are non-negative real numbers, but it is
sufficient to form the intersection using just a sequence of real
numbers $v < w$ which converges to $w$; similarly, in the last two
descriptions, it is sufficient to consider a sequence of rational
numbers $\alpha < w$ which converges to $w$.}
        \medskip
Proof: The second equality of these different descriptions of
$\Cal B_w$ is straight forward, the last one follows immediately
from Lemma 10. Let us show that $l(\bold p_w) =
\bigcap\nolimits_{v<w} \Cal Q_{v}.$ First, assume that $M$ belongs
to the intersection. Let $N$ be in $\bold p_w$. We want to  show
that $\Hom(M,N) = 0.$ There is a rational $\alpha$ with $0 <
\alpha < w$ such that $N$ belongs to $\bold p_\alpha$. Since $M$
is generated by $\bold t_\alpha$ and $\Hom(\bold t_\alpha,\bold
p_\alpha) = 0$, it follows that $\Hom(M,N) = 0$. Conversely,
assume that  $M$ belongs to $l(\bold p_w)$. Take a rational
$\alpha$ with $0 < \alpha < w$. We want to show that $M$ is
generated by $\bold t_\alpha$. Choose $\beta$ rational with
$\alpha <\beta < w$. Then $\bold t_\beta \subset \bold p_w$, thus
$l(\bold p_w) \subseteq l(\bold t_\beta) = g(G_\beta)$. And
$g(G_\beta) \subset g(\bold t_\alpha)$, according to Lemma 10.
This shows that $M$ is generated by $\bold t_\alpha$.
        \medskip
Proof of Theorem 6. For the second assertion, we only note that
$\Cal M(w) \subseteq \Cal C_w$ and that $\Cal M(w') \subseteq \Cal
B_{w'} \subseteq \Cal Q_w,$ since $w < w'$.

For the first assertion, let  $M$ be any indecomposable module
which does not belong to $\bold p_0$ or $\bold q_\infty$. Since
$M$ is indecomposable and does not belong to $\bold q_\infty$, we
have $\Hom(\bold q_\infty,M) = 0.$ Let $w$ be the infimum of all
$\alpha\in \Bbb Q_0^\infty$ such that $\Hom(\bold q_\alpha,M) =
0.$ Since $\bold q_w = \bigcup_{w < \alpha} \bold q_\alpha$, it
follows that $\Hom(\bold q_w,M) = 0,$ thus $M$ belongs to $\Cal
C_w$. It remains to be shown that $M$ also belongs to $\Cal B_w$. For $w = 0,$ this follows immediately from our assumption
that $M$ is indecomposable and does not belong to $\bold p_0$.
Thus, let $w > 0.$ We have to show that $M$ belongs to $\Cal
Q_\alpha$ for any rational number $\alpha$ with $0 < \alpha < w$.
Take such an $\alpha$ and assume that $M$ does not belong to $\Cal
Q_\alpha$. Since $(\Cal C_\alpha,\Cal Q_\alpha)$ is a split
torsion pair and $M$ is indecomposable, we conclude that $M$
belongs to $\Cal C_\alpha$. Thus $\Hom(\bold q_\alpha,M) = 0$. But
by the definition of $w$ this implies that $w \le \alpha$, a
contradiction.
        \bigskip
We add a further property of $\Cal M(w)$ which is quite useful to
know:
        \medskip
{\bf Proposition 9.} {\it The subcategories $\Cal C_w, \Cal B_w$
and $\Cal M(w)$ are closed under products and direct limits.}
        \medskip
We only have to consider the first two subcategories. Now $\Cal
C_w$ is the torsionfree class of a torsion pair, thus closed under
products. Also, since $\Cal C_w = r(\bold q_w)$, and $\bold q_w$
consists of finitely generated modules, we see that $\Cal C_w$ is
closed under direct limits.

Consider now $\Cal B_w$. All the subcategories $\Cal Q_w =
g(\bold q_w)$ are closed under direct limits, thus the same is
true for $\Cal B_w.$ It remains to be seen that $\Cal B_w$ is
closed under products.
 Assume that there are given modules $M_i \in \Cal B_w$
and let $M = \prod_{i\in I} M_i$.

Consider first the case $w > 0.$ Choose some $0 < \beta < w$ in
$\Bbb Q$. Then all the modules $M_i$ are generated by $G_\beta$,
thus there exist epimorphisms $\bigoplus_{I(i)} G_\beta \to M_i$
for some index set $I(i)$. But $\Add G_\beta$ is closed under
products, thus $\prod_{I(i)} G_\beta$ maps onto $\bigoplus_{I(i)}
G_\beta$ and therefore $\prod_i\prod_{I(i)} G_\beta$ maps onto
$\prod_i M(i)$. Again using that $\Add G_\beta$ is closed under
products, we see that $M$ is generated by $G_\beta$.

We proceed quite similarly for $w= 0.$ Write ${}_\Lambda \Lambda =
P\oplus P',$ where $P$ belongs to $\bold p_0$ and $P'$ to $\bold
t_0$ (thus $P = {}_\Lambda \Lambda_0$). Since $M_i$ does not split
off any indecomposable module from $\bold p_0$, any map $P \to
M_i$ can be factored through a module in $\add \tau^{-t}P$, for
any $t\in \Bbb N_0$. It follows that for any $t\in \Bbb N_0$, all
the modules $M_i$ are generated by $P_t = \tau^{-t} P \oplus P'.$
This module $P_t$ is a finite dimensional module, thus the
products of copies of $P_t$ are direct sums of copies of $P_t$.
This shows that $M$ itself is generated by $P_t$. Now consider an
indecomposable module $N$ from $\bold p_0.$ There is $t\in \Bbb
N_0$ with $\Hom(P_t,N) = 0$ and this implies that $\Hom(M,N) = 0.$
This shows that $M$ cannot split off a copy of $N$, as we had to
show.
        \medskip
{\bf Remark.} Note that in contrast to $\Cal B_w$, the
subcategories $\Cal Q_w$ are {\bf not} closed under products.
        \medskip
{\bf Examples.} For $\alpha\in \Bbb Q_0^\infty$, examples of
modules in $\Cal M(\alpha)$ have been mentioned above. Let us now
consider the case of an arbitrary $w\in \Bbb R^+$. In case $w$ is
not rational, $\Cal M(w)$ cannot contain any non-zero module of
finite length. We are going to provide two recipes for
constructing non-zero modules in $\Cal M(w).$

{\bf The first construction:} {\it Let $\alpha_1 > \alpha_2 >
\dots$ be a sequence of rational numbers converging to $w$ and
choose modules $M_i\in \add\bold t_{\alpha_i}.$ Then $\prod_i M_i
/ \bigoplus_i M_i$ belongs to $\Cal M(w)$.} Proof: Let $M =
\prod_i M_i$, and $M' = \bigoplus_i M_i$. Let us show that the
maximal submodule of $M$ which belongs to $\Cal Q_w$ is $M'$. On
the one hand, $M'$ is generated by $\bold q_w$. On the other hand,
given an indecomposable module $N \in \bold q_w$ and a non-zero
map $f\:N \to M$, then $N$ belongs to some $\bold t_\beta$ with
$\beta > w$. Since the sequence $(\alpha_i)_i$ converges to $w$,
there is some natural number $n$ with $\alpha_i < \beta$ for all
$i>n$. Thus the image of $f$ is contained in $\prod_{i \le n} M_i
\subseteq M'.$ Since $(\Cal C_w,\Cal Q_w)$ is a torsion pair, it
follows that $M/M'$ belongs to $\Cal C_w.$ In addition, we have to
show that $M/M'$ belongs to $\Cal B_w$. Since $\Cal B_w$ is
closed under products, $M$ belongs to $\Cal B_w$. But $\Cal B_w = l(\bold p_w)$ is closed under factor modules, thus with
$M$ also $M/M'$ belongs to $\Cal B_w.$

{\bf The second construction:} {\it Let $\alpha_1 < \alpha_2 <
\dots$ be a sequence of rational numbers converging to $w$ and
choose modules $M_i\in \add\bold t_{\alpha_i}$ with inclusions
$M_1 \subseteq M_2 \subseteq \dots$. Then the direct limit
$\lim\limits_\to M_i$  belongs to $\Cal M(w)$.} Proof: All the
modules $M_i$ belong to $\Cal C_w$ and $\Cal C_w$ is closed under
direct limits, therefore $M = \lim\limits_\to M_i$ belongs to
$\Cal C_w.$ Consider a rational number $\alpha$ where $0 < \alpha
< w$. There exists $i$ with $\alpha < \alpha_i$. Then $M_j$
belongs to $\add \bold q_\alpha$ for all $j \ge i$. This shows
that $M$ is generated by $\bold q_\alpha,$ and therefore $M$
belongs to $\Cal Q_\alpha$. As a consequence, $M$ belongs to $\Cal B_w$.

        \bigskip\bigskip
%======================================================================
{\bf References.}
        \medskip
\item{[AB]} Auslander, M., Buchweitz, R.O.: Maximal Cohen-Macaulay approximations.
  M\'emoire Soc\. Math\. France.  38 (1989), 5-37.
\item{[AC]} Angeleri-H\"ugel, L., Coelho, F\.U\.: Infinitely generated
 tilting modules of finite projective dimension. Forum Mathematicum. (To appear).
\item{[ATT]} Angeleri-H\"ugel, L., Tonolo, A., Trlifaj, J.: Tilting
 preenvelopes and cotilting precovers. Algebras and Representation Theory
 (To appear).
\item{[BS]} Buan, A. Solberg \O.: Limits of pure-injective
cotilting modules. (To appear).
\item{[CF]} Colby, R.R., Fuller, K.R.: Tilting, cotilting and serially tilted rings.
  Comm. Algebra 18 (1990), 1585-1615.
\item{[CET]} Colpi, R., D'Este, G., Tonolo, A.: Quasi-tilting modules and
counter equivalences. J.Algebra 191 (1997), 461-494. Corrigendum
in J.Algebra 206 (1998), 370.
\item{[CT]} Colpi, R., Trlifaj, J.: Tilting modules and tilting torsion theories.
 J.Algebra 178 (1995), 614-634.
\item{[DR]} Dlab, V., Ringel, C.M.: Indecomposable representations of
 graphs and algebras. Memoirs AMS 173 (1976).
\item{[F]} Fuchs, L, Infinite abelian groups, Vol. II, Academic
Press, New York 1973.
\item{[BK]} Buan, A.: Krause, H.: Cotilting modules over tame
hereditary algebras. (To appear)
\item{[HR]} Happel, D., Reiten, I.: An introduction to quasitilted
algebras, An. St.Univ.Ovidius Constantza, Vol. 4(1996) f.2,
137-149.
\item{[HRS]} Happel, D., Reiten, I., Smal\o, S.O., Tilting in abelian
categories and quasitilted algebras, Memoirs AMS, vol. 575 (1996)
\item{[K]} H. Krause: Generic modules over artin algebras.
  Proc. London Math. Soc.76 (1998), 276-306.
\item{[L]} Lenzing, H.: Generic modules over tubular algebras. In:
Advances in algebra and model theory. (G\"obel, ed.) Gordon \&
Breach (1997) 375-385.
\item{[LM]} Lenzing, H., Meltzer, H.: Tilting sheaves and concealed-canonical
algebras. ICRA VII Cocoyoc (Mexico) 1994. CMS Conf\. Proc\. 18
(1996), 455-473.
\item{[LP]} Lenzing, H., de la Pe\~na, J.: Concealed-canonical algebras and separating
 tubular families. Proc. London Math.Soc. 78(1999)513-540.
\item{[RS]} Reiten, I., Skowro\~nski, A.: Sincere stable tubes, J.
Alg. 232(2000) 64-75.
\item{[R1]} Ringel, C.M.: Infinite dimensional representations of finite
 dimensional hereditary algebras. Symposia Math. XXIII. Istituto Naz\. Alta Mat\. (1979),
 321-412.
\item{[R2]} Ringel, C.M.: Tame algebras and integral quadratic forms.
 Springer LNM 1099 (1984).
\item{[R3]} Ringel, C.M.: Representation theory of finite
dimensional algebras, Durham lectures 1985, London Math. Soc.
Lecture Note Series 116 (1986) 7-79.
\item{[R4]} Ringel, C.M.: Canonical algebras (with an appendix by W.W. Crawley Boevey).
  In: Topics in Algebra. Banach Center Publications. 26 Part 1 (1990) 407-432.
\item{[R5]} Ringel, C.M.: The Ziegler spectrum of a tame
hereditary algebras, Coll. Math. 76 (1998) 105-115.
\item{[R6]} Ringel, C.M.: A Construction of Endofinite Modules. in: Advances in Algebra
 and Model Theory. Gordon-Breach. London (1997), 387-399.
\item{[R7]} Ringel, C.M.: Tame algebras are wild. Algebra Colloq. 6 (1999),
 473-480.
\item{[R8]} Ringel, C.M.: Infinite length modules. Some examples as introduction. In:
 Infinite Length Modules (Krause, Ringel ed.) Birkh\"auser (2000), 1-73.

        \bigskip
\noindent

%===================================
\baselineskip=9pt \leftskip=1truecm
\parindent=0truecm
\rmk

I.Reiten, Department of Mathematical Sciences, Norwegian
University of Science and Technology,
\par 7490 Trondheim, Norway\par
E-mail address: {\ttk reiten\@math.ntnu.no}
        \medskip
C.M.Ringel, Fakult\"at f\"ur Mathematik, Universit\"at
Bielefeld,\par POBox 100\,131, \ D-33\,501 Bielefeld, Germany\par
E-mail address: {\ttk ringel\@mathematik.uni-bielefeld.de}

\bye